\input ./style/arxiv-general.cfg
\documentclass[aop,MSNbibl,dvips]{arximspdf}
\makeatletter
   \@ifpackageloaded{graphicx}{}{\usepackage{graphicx}}
\makeatother

\doi{10.1214/14-AOP945} 
\volume{43}
\issue{5}
\pubyear{2015}
\firstpage{2611}
\lastpage{2646}
\docsubty{FLA}

\makeatletter
\newcommand{\bull}{\bullet}
\newtheorem{prop}{Proposition}
\newproclaim{defn}[prop]{Definition}
\newtheorem{lm}[prop]{Lemma}
\newproclaim{ex}[prop]{Example}
\newtheorem{teo}[prop]{Theorem}
\newtheorem{cor}[prop]{Corollary}
\newproclaim{rem}[prop]{Remark}
\newcommand{\hatF}{\widehat{F}}
\newcommand{\hatM}{\widehat{M}}
\newcommand{\tilM}{\widetilde{M}}
\newcommand{\hatp}{\widehat{p}}
\newcommand{\giv}{|}
\makeatother

\begin{document}
\begin{frontmatter}

\title{Regenerative tree growth: Markovian embedding of~fragmenters,
bifurcators, and bead splitting~processes\thanksref{T1}}
\runtitle{Regenerative tree growth}

\begin{aug}
\author[A]{\fnms{Jim}~\snm{Pitman}\ead[label=e1]{pitman@stat.berkeley.edu}}
\and
\author[B]{\fnms{Matthias}~\snm{Winkel}\corref{}\ead[label=e2]{winkel@stats.ox.ac.uk}}
\runauthor{J. Pitman and M. Winkel}
\affiliation{University of California Berkeley and University of Oxford}
\address[A]{Department of Statistics\\
University of California, Berkeley\\
Berkeley, California 94720\\
USA\\
\printead{e1}}
\address[B]{Department of Statistics\\
University of Oxford\\
1 South Parks Road\\
Oxford OX1 3TG\\
United Kingdom\\
\printead{e2}}
\end{aug}
\thankstext{T1}{Supported by the National Science Foundation Awards
0405779 and 0806118.}

\received{\smonth{4} \syear{2013}}
\revised{\smonth{5} \syear{2014}}

%
\begin{abstract}
Some, but not all processes of the form $M_t=\exp(-\xi_t)$ for a
pure-jump subordinator $\xi$ with
Laplace exponent $\Phi$ arise as residual mass processes of particle 1
(tagged particle) in Bertoin's partition-valued exchangeable
fragmentation processes. We
introduce the notion of a \emph{Markovian embedding} of $M=(M_t,t\ge
0)$ in a fragmentation process, and we show that for
each $\Phi$, there is a unique (in distribution) binary fragmentation
process in which $M$ has a Markovian embedding. The identification of
the Laplace exponent $\Phi^*$ of its tagged particle process $M^*$
gives rise to a symmetrisation operation $\Phi\mapsto\Phi^*$, which
we investigate in a general study of pairs $(M,M^*)$ that
coincide up to a random time and then evolve independently. We call $M$
a \emph{fragmenter} and $(M,M^*)$ a
\emph{bifurcator}.

For $\alpha>0$, we equip the interval $R_1=[0,\int_0^\infty
M_t^\alpha \,dt]$
with a purely atomic probability measure $\mu_1$, which captures the
jump sizes of $M$ suitably placed on $R_1$. 
We study binary tree growth
processes that in the $n$th step sample an atom (``bead'') from $\mu
_n$ and build $(R_{n+1},\mu_{n+1})$ by replacing the atom by a
rescaled independent copy of $(R_1,\mu_1)$ that we tie to the position
of the atom. We show that any such
\emph{bead splitting process} $((R_n,\mu_n),n\ge1)$ converges almost
surely to an $\alpha$-self-similar continuum random tree of Haas and
Miermont, in the Gromov--Hausdorff--Prohorov
sense. This generalises Aldous's line-breaking construction of the
Brownian continuum random tree.
\end{abstract}

%
\end{frontmatter}

\section{Introduction}\label{sec1}

We call a process $M:= (M_t, t \ge0)$ a \emph{multiplicative
subordinator}, or \emph{fragmenter} for short, if
\[
M_t = \exp( - \xi_t ),\qquad t \ge0,
\]
for some subordinator $(\xi_t, t \ge0)$. As shown by Pitman \cite
{Pit-99} and Bertoin \cite{Ber-hom}, such
processes arise naturally in the theory of continuous-time processes of
coagulation and
fragmentation.
The process $(1 - M_t, t \ge0 )$ is the random cumulative distribution
function of
a random discrete probability measure on $(0,\infty)$. These random
measures have been
studied in the theory of Bayesian nonparametric statistics \cite{Dok-74,DoJ-04,Jam-06}, not just for
subordinators $\xi$, but more generally for increasing processes with
independent
increments which are not necessarily stationary.

We will use terminology based on the fragmentation interpretation of
$M_t$ as the residual mass of a block containing a particle at time
$t$. Bertoin \cite{Ber-hom} showed that the mass containing particle 1
in an exchangeable homogeneous fragmentation process is a fragmenter.
Let us recall the definition of a homogeneous fragmentation process
(HFP). We denote by $\mathcal{P}=\mathcal{P}_\mathbb{N}$ the set of
partitions of $\mathbb{N}$. An exchangeable HFP is a Markov process
$\Pi=(\Pi(t),t\ge0)$ in $\mathcal{P}$ such that:
\begin{itemize}
\item given $\Pi(t)=\{B_i,i\ge1\}\in\mathcal{P}$, the partition
$\Pi(t+s)$ is distributed as the collection of blocks of $B_i\cap\Pi
^{(i)}(s)$, $i\ge1$, for a family $\Pi^{(i)}$, $i\ge1$, of
independent copies of $\Pi$, and
\item the distribution of $\Pi$ is exchangeable, that is, invariant
under all finite permutations of $\mathbb{N}$.
\end{itemize}
It is a well-known consequence of de Finetti's theorem that
exchangeable partitions have asymptotic frequencies, so $|\Pi
_i(t)|=\lim_{n\rightarrow\infty}n^{-1}\#(\Pi_i(t)\cap\{1,\ldots,n\})$
exists almost surely (which we abbreviate \emph{a.s.}), in fact
jointly for all $i\ge1$ and~$t\ge0$; see~\cite{Ber-hom}. Referring to an asymptotic frequency as \emph{mass},
$M^*_t:=|\Pi_1(t)|$ is the residual mass of the block $\Pi_1(t)$,
which contains particle 1. Bertoin showed that for every exchangeable
HFP $\Pi$, this process $M^*$ is a fragmenter. We will call $M^*=|\Pi
_1|:=(|\Pi_1(t)|,t\ge0)$ the \emph{canonical fragmenter} of $\Pi$.
As Haas \cite{Haa-03} demonstrated, there are exchangeable HFPs with
different distributions whose canonical fragmenters have the same
distribution. On the other hand, in the subclass of binary models,
where every infinitesimal split is into two parts (see Section~\ref{stexhom}), we show this cannot happen.

The starting point for this paper is the observation that not all
fragmenters arise as canonical fragmenters in an exchangeable HFP.
Furthermore, we have encountered a number of natural nonexchangeable
models \cite{For-05,CFW,PW09,CW,PRW}, in which masses of blocks can be
defined as asymptotic frequencies, and the mass containing particle 1
is also a fragmenter. Via embedding of such residual mass processes
into an exchangeable model or via limit considerations, we have found
associated exchangeable models in all those examples. Our main result,
Theorem~\ref{embed}, shows that for a suitable notion of
``embedding,'' these examples generalise to a remarkably simple picture:
\begin{longlist}[1.]
\item[1.] Every pure-jump fragmenter $M$ can be embedded in an
exchangeable \mbox{binary} HFP $\Pi$.
\item[2.] The distribution of the HFP $\Pi$ is uniquely determined
by that of the fragmenter $M$.
\item[3.] The canonical fragmenter $M^*$ of $\Pi$ is a symmetrised
version of $M$ defined in Section~\ref{sizebiased}.
\end{longlist}
To\vspace*{1pt} prepare this result, Section~\ref{bifurcators} offers a systematic
study of pairs of fragmenters $(M,\hatM)$ that coincide up to a random
time, after which $M$ and $\hatM$ evolve independently. We call such
pairs \emph{bifurcators} and give several equivalent characterisations
(Propositions~\ref{propsplit1} and~\ref{propsplit2}), which are of
interest in their own right. Examples of bifurcators include pairs
$((|A^1_t|,|A^2_t|),t\ge0)$ of residual masses of the blocks $A^1_t$
and $A^2_t$ containing particles 1 and 2, respectively, in an
exchangeable HFP; see also Proposition~\ref{symmbif}. Specifically,
note that after the random time when $1$ and $2$ separate, the
evolution of $A^1_t$ and $A^2_t$ is independent, since disjoint blocks
evolve independently in a HFP. In Section~\ref{sizebiased} we focus
more generally on bifurcators for which size-biased switching describes
the separation time. This induces an idempotent transformation from
distributions of $M$ to $M^*$, which we call \emph{symmetrisation}.
This transformation is the key to finding the exchangeable binary HFP
\mbox{associated} with~$M$. In Section~\ref{stexhom} we recall-known facts
about exchangeable HFPs~$\Pi$. We introduce the notion of a \emph{Markovian embedding in $\Pi$} in Section~\ref{fragemb} and show
the existence of an embedding for $M$. We postpone the proof of
uniqueness of $\Pi$ to Section~\ref{genline}. In Section~\ref{sectmass} we study the three-way mass split into the parts before and
after the separation time, and in Section~\ref{sectlength} associated
lengths induced by a bifurcator.

To complete the proof of Theorem~\ref{embed}, we use Haas and
Miermont's \cite{HM04} $\alpha$-self-similar continuum random trees
({CRTs}), which are certain random rooted compact metric space
trees $(\mathcal{T},d,0)$ equipped with a probability measure $\mu$.
Specifically, a random element $\Sigma_1^*\in\mathcal{T}$ with
distribution $\mu$ yields a path $\mathcal{R}_1^*=[\![0,\Sigma_1^*]\!]$
in $\mathcal{T}$. For $(\mathcal{T},\mu)$ associated with $\Pi$
(and some $\alpha>0$), it is well known that the process of $\mu
$-masses in subtrees above points in $\mathcal{R}_1^*$ is\vspace*{1pt} related to a
copy of $M^*$ by a certain $\alpha$-self-similar time change.
Furthermore, $\mathcal{R}_n^*=\bigcup_{j=1}^n[\![0,\Sigma_j^*]\!]$
increases to $\mathcal{T}$ for a sample $\Sigma_n^*$, $n\ge1$, from
$\mu$. See Section~\ref{stexhom} for details.

We project $\mu$ onto $\mathcal{R}_n^*$ to equip $\mathcal{R}_n^*$
with a random discrete distribution $\mu_n^*$. In particular,
$(\mathcal{R}_1^*,\mu_1^*)$ is a \emph{string of beads}, that is, an
interval equipped with a purely atomic measure, and $(\mathcal
{R}_1^*,\mu_1^*)$ can easily be completely expressed in terms of
$M^*$; see Section~\ref{sectbsp}. Note that $\mathcal{R}_{n+1}^*$ is
a tree with one more branch than $\mathcal{R}_n^*$, and $\mu_n^*$ is
the projection of $\mu_{n+1}^*$ onto $\mathcal{R}_n^*$. Since $\Sigma
_{n+1}^*$ is selected according to $\mu$, we have an instance of the
following general notion of a \emph{bead splitting process}:
\begin{itemize}
\item Let $(R_1,\mu_1)$ be a string of beads.
\item Given that $R_n$ has been defined with a purely atomic
probability measure $\mu_n$, pick an atom (``bead'') $J_n$ according
to $\mu_n$. Given $\mu_n(\{J_n\})=m$, remove the atom from $\mu_n$,
split it into smaller atoms, tie to $J_n$ a string with
these beads of total mass $m$ to form $(R_{n+1},\mu_{n+1})$.
\end{itemize}

%
\begin{figure}

\includegraphics{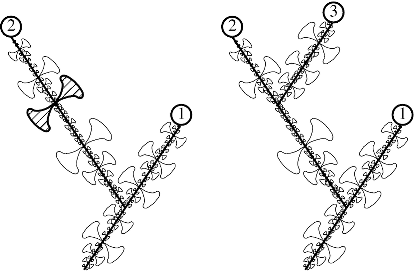}

\caption{Binary trees equipped with strings of beads; the right tree
was obtained from the left tree by performing the above operation on
the shaded bead of the left tree.}\label{fig1}
\end{figure}

See Figure~\ref{fig1} for an illustration. A similar bead splitting
process, but with different bead selection rules, was obtained for the
alpha--theta model of \cite{PW09} by exploiting properties of the
Chinese restaurant process. The main developments in Section~\ref{genline} culminate in Theorem~\ref{teogrowth}:
\begin{longlist}[3.]
\item[4.] We give an autonomous description meaningful outside a CRT
for the evolution of the bead-splitting process $((\mathcal{R}_n^*,\mu
_n^*),n\ge1)$, which is associated with the canonical
fragmenter.
\item[5.] We generalise this description to start from a string of
beads associated with an arbitrary strictly decreasing pure-jump
fragmenter $M$ rather than a canonical fragmenter $M^*$. At each growth
step for $(R_{n+1},\mu_{n+1})$ we tie to $J_n$ an independent rescaled
copy of $(R_1,\mu_1)$.
\item[6.] We show that this bead splitting process $((R_n,\mu
_n),n\ge1)$ converges almost surely for the
Gromov--Hausdorff--Prohorov metric to a CRT $(\mathcal{T},\mu)$. This
CRT is associated with a HFP whose canonical fragmenter is $M^*$ as
identified earlier.
\end{longlist}
As tools, we develop a general spinal decomposition of exchangeable
HFPs along a Markovian path (Lemma~\ref{spindecA}) and show a CRT
convergence result for bead splitting processes based on any Markovian
path (Lemma~\ref{markpathconv}), which we then also use to complete
the proof of Theorem~\ref{embed}. The embedding for the existence part
of the proof of Theorem~\ref{embed} is not carried out in a CRT, but
directly in an exchangeable HFP. A uniqueness proof entirely in the
framework of HFP should be possible, but the construction is harder to
formulate, and the compactness of CRTs would not be directly available.
The transition kernel from $(R_n,\mu_n)$ to $(R_{n+1},\mu_{n+1})$ is
simple for all fragmenters and gives an inductive description of the
distribution of $(R_n,\mu_n)$ for every $n\ge1$. Section~\ref{bcrt}
provides direct descriptions of the distribution of $(\mathcal
{R}_n^*,\mu_n^*)$ in the special case of the Brownian CRT, exploiting
relations to Aldous's line-breaking construction \cite{Ald-crt1},
Brownian path transformations \cite{AMP-04a} and Poisson--Dirichlet
distributions \cite{DGM,csp}.

The CRT convergence result of bead splitting processes here can be
complemented by scaling limit results of discrete tree shapes $T_n$ of
$R_n$ and/or their reduced subtrees $T_{n,k}$ spanned by the first $k$
leaves. Specifically, we applied methods of Haas and Miermont \cite
{HM10} in \cite{PRW} to obtain convergence in distribution for trees
like $T_n$, $n\ge1$, suitably rescaled, to a limiting CRT.

We note in Proposition~\ref{regcomp} that the numbers of leaves of
$\mathcal{R}_{n+1}$ in subtrees of the spine from the root to leaf 1
form a strongly sampling consistent regenerative composition structure
$\mathcal{C}_n$, $n\ge1$, in the sense of Gnedin and Pitman~\cite{gp}, also~\cite{PW09}, Section~2.1. Gnedin, Pitman, and Yor \cite
{GPYalpha} studied the number of blocks $\#\mathcal{C}_n$ and showed
$\#\mathcal{C}_n/n^\alpha\rightarrow\int_0^\infty M_s^\alpha \,ds$
a.s., under a
regular variation condition. We exploited this in \cite{HMPW},
Proposition~7, for exchangeably labelled trees to see that reduced trees
$T_{n,k}$ converge almost surely to $R_k$ when rescaling all
edge lengths by $n^\alpha$. This result can be generalised to the
present setup, which includes nonexchangeable cases. See also \cite
{PW09}, Proposition~14, for the alpha--theta model, which adds projected
uniform measures that converge to the limiting strings of beads. The
bead splitting process we identified for the alpha--theta model develops
by size-biased branching only for $\theta=\alpha$. For the other
cases, we found different bead selection rules in connection with
ordered Chinese restaurant processes.

\section{Fragmenters and their embedding in fragmentation processes}\label{sect4}

This section studies a natural class of models for tracking two
residual mass processes $(M_t,t\ge0)$ and $(\hatM_t,t\ge0)$ that we
can think of
as parts in a fragmentation process. It is instructive and indeed a
natural approach to fragmentation processes to first consider these
models in their own right, not as parts in a fragmentation process. A
systematic study does not appear to be available in the literature.
This is provided in Section~\ref{bifurcators}, before Section~\ref
{sizebiased} focusses on the important special case of size-biased
branching. In Section~\ref{stexhom} we recall from the literature the
concepts of a self-similar fragmentation process and associated CRTs.
This enables us in Section~\ref{fragemb} to apply results from
Sections~\ref{bifurcators} and~\ref{sizebiased} to formalise and
establish points 1 and 3 from the \hyperref[sec1]{Introduction}. Sections~\ref
{sectmass} and~\ref{sectlength} study the three-way mass split and
associated lengths before and after the separation time of a bifurcator
under size-biased branching.

\subsection{Fragmenters, switching transformations, and bifurcators}\label{bifurcators}

It will be assumed throughout this section that all fragmenters
$M_t=\exp(-\xi_t)$, $t\ge0$,
are derived from subordinators $\xi$ with zero drift and no killing.
Furthermore, for most of our discussion, we will also assume an
absolutely continuous L\'evy density. This is just for convenience of
presentation. We discuss general L\'evy measures at the end of this
section. 
The L\'evy--It\^o representation of $\xi$ is then
\[
\xi_t = \sum_{0 < s \le t } \Delta
\xi_s, \qquad t \ge0, %
\]
where $\{ (s, \Delta\xi_s )\dvtx  s >0, \Delta\xi_s > 0 \}$ is the set
of points of
a Poisson random measure on $(0,\infty) \times(0,\infty)$ with
intensity measure
$ ds\, \lambda(x)\,dx$ where $\lambda$ is the L\'evy density of the
subordinator, so
%
\begin{equation}\label{LK}
\quad\mathbb{E}\bigl( e^{- \rho\xi_t } \bigr) = e^{- t \Phi(\rho)
}\qquad\mbox{where } \Phi(\rho) = \int_0^\infty\bigl(1 -
e^{- \rho x } \bigr) \lambda( x ) \,dx, \rho\ge0,
\end{equation}
is the Laplace exponent.
Let $F_s:= \exp( - \Delta\xi_s )$. Then the corresponding
formulae for the fragmenter $M$ are
\[
M_t = \prod_{0 < s \le t } F_s,
\qquad t \ge0, %
\]
where $\{ (s, F_s )\dvtx  s>0, F_s < 1 \}$ is the set of points of
a Poisson random measure on $(0,\infty) \times(0,1)$ with intensity measure
$ ds\, uf(u)\,du$ on $(0,\infty)\times(0,1)$, where $uf(u)\,du$ is the
push-forward of
$\lambda(x)\,dx$ via the transformation $u = e^{-x}$. So for all
nonnegative Borel functions $g$,
%
\begin{equation}
\label{X} \int_0^\infty g\bigl(e^{-x}
\bigr)\lambda(x)\,dx=\int_0^1 g(u)uf(u)\,du.
\end{equation}
%
We introduce the size-biasing factor $u$ in the definition (\ref{X}) of
$f(u)$ to simplify applications to fragmenters associated with (binary)
homogeneous fragmentations~\cite{ber-book}, which we define more
formally in Section~\ref{stexhom} and explain briefly in the next
paragraph. We call $f$ the \emph{splitting density of the fragmenter},
which is related to $\lambda$
by
%
\begin{equation}
\label{fla} f(u) = u^{-2} \lambda( - \log u ), \qquad0 < u < 1.
\end{equation}
By (\ref{X}) for $g(u)=1-u^\rho$, the Laplace exponent of $\xi$ is then
%
\begin{equation}
\label{lapex} \Phi(\rho) = \int_0^1 \bigl( 1 -
u ^\rho\bigr) u f(u) \,du, \qquad\rho\ge0.
\end{equation}
Note that $f$ is subject to the integrability condition that $\Phi
(\rho) < \infty$
for some (hence all) $\rho> 0$, that is,
%
\begin{equation}
\Phi(1) = \int_0^1 u ( 1 - u)f(u) \,du < \infty.
\end{equation}
The L\'evy--Khintchine formula (\ref{LK}) now provides a Mellin transform
for the fragmenter,
\[
\mathbb{E}\bigl( M_t^\rho\bigr) = e^{- t \Phi(\rho) }, \qquad
t \ge0. %
\]
%
If $(M_t, t \ge0 )$ is the mass of a randomly
tagged fragment in a binary homogeneous fragmentation process with a dislocation
measure $\nu_{\mathrm{ranked}}$ concentrated on decreasing nonnegative
sequences $\mathbf{s}=(s_1, s_2, \ldots)$
with $s_1 + s_2 = 1$, then $M$ admits the above descriptions,
assuming the existence of a density $\sum_{i=1}^2 \nu_{\mathrm{ranked}}(s_i \in du ) =f(u) \,du$. The size-biasing factor $u$
then arises because $f$ is necessarily \emph{symmetric}, meaning $f(u)= f(1-u)$,
and given a mass split $\mathbf{s}=(u,1-u)$ with $u > 1-u$ from $\nu
_{\mathrm{ranked}}(d\mathbf{s})$, the randomly tagged fragment will be found in
component $u$ with probability $u$ and in component $1-u$ with
probability $1-u$. This is mapping a Poisson point process of mass
splits $\nu_{\mathrm{ranked}}(d\mathbf{s})$ to a Poisson point process with
intensity $uf(u)\,du$.

To further study fragmenters embedded in a homogeneous fragmentation
process, we consider the following \emph{switching transformation} of
one fragmenter $M$ into another fragmenter $\hatM$.
Let $p$ be a nonnegative measurable function from $(0,1)$ to $[0,1]$.
If $M$ is a fragmenter with L\'evy--It\^o representation
$M_t = \prod_{0 < s \le t } F_s$, consider the
process
%
\begin{equation}
\label{itohat} \hatM_t:= \prod_{0 < s \le t }
\hatF_s,
\end{equation}
where conditionally given $M$, the factors
$\hatF_s$ are defined by $\hatF_s = F_s$ with probability
$1-p(F_s)$, and
$\hatF_s = 1 - F_s$ with probability $p(F_s)$.
Here the construction of the point process $((s,\hatF_s), s > 0 )$ from
the point process $((s,F_s), s > 0 )$ with intensity $ds\, u f(u) \,du$
is made rigorous in the usual way by some arbitrary indexing
of these points by positive integers, and making independent choices
for each of the
countable number of $F_s$ with $F_s < 1$.
Here and below it is always
assumed that all processes are defined on a rich enough probability space
to admit all necessary auxiliary\vspace*{1pt} randomizations, as are involved in
passing from
$((s,F_s), s > 0 )$ to $((s,\hatF_s), s > 0 )$.
Standard transformation results for Poisson\vspace*{1pt} point
processes imply that $((s,\hatF_s), s > 0 )$ is a Poisson point
process with
intensity $ds\, u \widehat{f}(u) \,du$ for
$\widehat{f}$ determined by the formula
%
\begin{equation}
\label{hatf} \quad u \widehat{f}( u ) = \bigl(1 - p(u)\bigr) u f(u) + p(1-u) (1-u)
f(1-u), \qquad0 < u < 1.
\end{equation}
In particular, provided $\int_0^1 u (1-u) \widehat{f}(u)\,du < \infty
$, that is,
%
\begin{equation}
\label{pcond} \int_0^1 u^2 p(u)
f(u) \,du < \infty,
\end{equation}
the function $\widehat{f}$ serves as a splitting density,
and (\ref{itohat}) is the L\'evy--It\^o representation of
a fragmenter $\hatM$ with splitting density $\widehat{f}$.
Call $\hatM$ the \emph{fragmenter derived from~$M$ by switching
according to $p$.}

The following proposition provides a summary:

%
\begin{prop}
\label{propsplit}
If $M$ is a fragmenter with splitting density $f$, and $p$ is subject
to (\ref{pcond}),
then $\hatM$ derived from
$M$ by switching according to $p$ is a fragmenter with splitting density
$\widehat{f}$ as in (\ref{hatf}).
Moreover, $M$ is then derived from $\hatM$ by switching according to
$\widehat{p}$,
where
%
\begin{equation}
\label{hatp} \widehat{p}(u) = \frac{ p(1-u) (1-u) f(1-u)} { u \widehat
{f} (u) }, \qquad0 < u < 1.
\end{equation}
\end{prop}

This generalises \cite{PW09}, Lemma 19(b), which treats the case of
size-biased switching probabilities $p(u)=1-u$.

Observe that two fragmenters $M$ and $\hatM$ as above are by
construction such
that $M_t = \hatM_t$ for $ 0 \le t < \tau$ where
\[
\tau:= \inf\{ s\dvtx  F_s \ne\widehat{F}_s \}
\]
is the time of the first switch.
It is clear from the Poisson construction that $\tau$ is exponential with
rate
%
\begin{equation}
\label{phiphihat} \phi:= \int_0^1 p(u) u f(u) \,du
= \int_0^1 \hatp(u) u \widehat{f}(u) \,du=:
\widehat{\phi}\in[0,\infty],
\end{equation}
where we will usually exclude the trivial cases $\phi=0$, that is,
$\tau=\infty$, and $\phi=\infty$, that is, $\tau=0$.
The conditional distribution of $\tau$ given $M$ is made
explicit by the formula
%
\begin{equation}
\label{taugivM} \mathbb{P}( \tau> t \giv M ) = \prod
_{0 < s \le t } \bigl( 1 - p(F_s) \bigr),
\end{equation}
where $F_s = M_{s}/M_{s-}$ and, by convention, $p(1)=0$.
Assuming further that $0<\phi< \infty$, so $\mathbb{P}(0 < \tau<
\infty) = 1$,
it is clear by construction that $M$ and $\hatM$ satisfy the
\emph{identification rule}
%
\begin{equation}
\label{branch} M_t = \hatM_t\qquad\mbox{for } 0 \le t <
\tau,
\end{equation}
hence $M_{\tau-} = \hatM_{\tau-}$, and the \emph{binary splitting condition}
that the decrement of each fragmenter at time $\tau$
equals the value of the other fragmenter at time~$\tau$,
%
\begin{equation}
\label{switch} M_{\tau_-} - M_{\tau} = \hatM_{\tau} \quad\mbox{and}\quad\hatM_{\tau_-} - \hatM_{\tau} = M_{\tau}.
\end{equation}
Call (\ref{branch}) and (\ref{switch}) together the \emph{binary
junction conditions}.
After time $\tau$ the random factors governing the evolution of $M$
and $\hatM$
are further coupled. We now modify this construction so that the two
fragmenters continue independently after time $\tau$:

%
\begin{defn}[(Bifurcator)]\label{bifdef}
We call a pair of
fragmenters $(M,\tilM)$ a \emph{bifurcator with switching time $\tau$}
if there are
a splitting density $f$ and a switching probability function $p$ so
that $\int_0^1p(u)uf(u)\,du<\infty$ and $(M,\tilM,\tau)$ has the
following joint distribution:
\begin{itemize}
\item$M$ is a fragmenter with splitting density $f$,
\item$\tau$ is the first switching time of an auxiliary fragmenter
$\hatM$ derived from $M$ by switching according to $p$,\vspace*{1pt}
\item$\tilM_t = M_t$, $t < \tau$, $\tilM_\tau= M_{\tau-}
- M_\tau$, and $(\tilM_{\tau+t}/\tilM_\tau,t \ge 0)$ is a
copy of $\hatM$
independent of \nolinebreak$(M,\tau)$.
\end{itemize}
\end{defn}

See Propositions~\ref{propsplit1} and~\ref{propsplit2} for
characterisations that may serve as alternative definitions.

Note that in our construction, $\tilM= \hatM$ on $[0,\tau]$, so the
binary junction conditions (\ref{branch}) and (\ref{switch}) hold
just as well for $\tilM$ as for $\hatM$. But after time $\tau$
the evolutions of $M$ and $\tilM$ are decoupled. Dually, $(M_{\tau
+t}/M_{\tau}, t \ge0 )$
is a copy of $M$ which is independent of $(\tilM,\tau)$.

Henceforth we will no longer be concerned with any $\hatM$ that is
further coupled with $M$ after $\tau$, and
we will instead use the generic notation $(M,\hatM)$ for a bifurcator.
Then for some splitting time $\tau$, whose joint law with $M$ is
determined by the switching
probability function $p$,
%
\begin{equation}
\hatM_t = M_t 1_{\{\tau>t\}} + ( M_{\tau-} -
M_{\tau}) \hatM'_{t-\tau} 1_{\{\tau\le t\}},
\end{equation}
where $\hatM' \stackrel{d}{=}\hatM$ with $\hatM'$ independent of
$(M,\tau)$.
Note the subtlety that $\hatM$ is determined pathwise by $M$ up to and
including the
splitting time $\tau$, but thereafter the jumps of
$\hatM$ and $M$ are decoupled: the distribution of how $\hatM$
evolves after time
$\tau$ is implicitly determined
by $M$ and $p$, but\vspace*{1pt} there is no pathwise coupling between $M$ and
$\hatM$ after time
$\tau$. Rather, $(M,\hatM)$ satisfies

%
\begin{defn}[(Asymmetric Markov branching property)]\label{ambp} We
say that $(M,\hatM)$ has the \emph{asymmetric Markov branching property}
relative to the splitting time $\tau$ if:
\begin{itemize}
\item conditionally given $\tau> t$ the process $((M_{t+v}/M_t,\hatM
_{t+v}/\hatM_t), v \ge0)$ is a
copy of $(M,\hatM)$, independent of $((M_s,\hatM_s), 0 \le s \le t)$;
\item conditionally\vspace*{1pt} given $\tau\le t$ the two processes $(M_{t+v}/M_t,
v \ge0)$ and $((\hatM_{t+v}/\break \hatM_t), v \ge0)$
are independent copies of $M$ and $\hatM$, respectively, independent of
$((M_s,\hatM_s),0\le s \le t)$.
\end{itemize}
\end{defn}

The following variation of Proposition~\ref{propsplit} follows easily
from standard facts about Poisson point processes, and the above definitions:

%
\begin{prop}\label{propsplit1}
\textup{(a)}
The joint distribution of the bifurcator $(M,\hatM)$ is uniquely
determined by the splitting density $f$ of $M$ and a switching
probability function $p$
subject to
%
\begin{equation}
0<\phi:= \int_0^1 p(u) u f(u) \,du < \infty
\end{equation}
or dually by the splitting density $\widehat{f}$ of $\hatM$ and the
dual switching probability
function $\hatp$, subject to $0<\widehat{\phi} < \infty$, as
specified in (\ref{hatf}), (\ref{hatp}) and (\ref{phiphihat}).
Furthermore, $\phi= \widehat{\phi}$ is the rate of the exponentially
distributed junction time $\tau$.

\textup{(b)}
A bifurcator $(M,\hatM)$ as in \textup{(a)} can also be constructed as
follows from five independent ingredients:\vspace*{1pt}
three fragmenters $M^0,M'$, and $\hatM'$ with splitting densities $(1
- p(u) ) f(u)$,
$f(u)$, and $\widehat{f}(u)$, respectively, an exponential time $\tau$
with rate~$\phi$, and
a random variable $U \in(0,1)$ with distribution
%
\begin{equation}
\label{vdens} \mathbb{P}( U \in du ) = \phi^{-1} p(u) u f(u) \,du,
\qquad0 < u < 1.
\end{equation}
Now define $(M,\hatM)$ by $M_t = \hatM_t = M^0_t$ for $t<\tau$, and let
%
\begin{equation}
\label{vfac} M_{\tau+ v } = U M^0_{\tau}
M_v^\prime\quad\mbox{and}\quad\hatM_{\tau+ v } =
(1-U ) M^0_{\tau} \hatM_v^\prime\qquad\mbox{for } v \ge0.
\end{equation}
%
%
\end{prop}

We assumed for ease of exposition that $M$ has a splitting density $f$.
However, the operation of switching according to $p$ and the notion of
an associated bifurcator are meaningful when we replace $uf(u)\,du$ by a
more general measure $\Lambda(du)$
satisfying $\int_{(0,1)}(1-u)\Lambda(du)<\infty$.
We generalise (\ref{hatf}) and (\ref{hatp}) to
%
\begin{eqnarray}\label{hatlambda}
\widehat{\Lambda}(du) &=& \bigl(1-p(u)\bigr)\Lambda(du)+p(1-u)
\overline{\Lambda}(du)\quad\mbox{and}
\nonumber\\[-8pt]\\[-8pt]
\widehat{p}(u)\widehat{\Lambda}(du)&=&p(1-u)\overline{\Lambda}(du),\nonumber
\end{eqnarray}
where $\overline{\Lambda}$ is the image measure of $\Lambda$ under
the switching operation $u\mapsto1-u$. Then (\ref{hatlambda}) defines
$\widehat{\Lambda}(du)$ as a measure satisfying
$\int_{(0,1)}(1-u)\widehat{\Lambda}(du)=2\int_{(0,1)}p(u)\Lambda
(du)<\infty$.
Also, $p(1-u)\overline{\Lambda}(du)$ is, by definition of $\widehat
{\Lambda}$, absolutely continuous with respect to $\widehat{\Lambda
}$ with density taking values in $[0,1]$. This identifies $\widehat
{p}(u)$ for $\widehat{\Lambda}$-a.e. $u\in(0,1)$. We define $\tau$
as the first switching time with distribution given in (\ref{taugivM}).
If $p(1/2)\Lambda(\{1/2\})>0$, we can have $\tau\neq\inf\{s\ge
0\dvtx  F_s\neq\widehat{F}_s\}$, and if furthermore $\Lambda$ is
finite, also
$\tau\neq\inf\{s\ge0\dvtx  M_s\neq\widehat{M}_s\}$ for a
bifurcator $(M,\widehat{M})$. A~more satisfactory way to include those
cases is to slightly re-model the switching transformation by marking
$(F_s,s\ge0)$ by a marking kernel $K$ from $(0,1)$ to $\{0,1\}$, where
$K(u,\{1\})=p(u)$ and $K(u,\{0\})=1-p(u)$, with associated marked point
process $((F_s,m_s),s\ge0)$ mapping to $\widehat
{F}_s=(1-m_s)F_s+m_s(1-F_s)$, and with $\tau=\inf\{s\ge0\dvtx
m_s=1\}$. In the following characterisation of bifurcators, it is more
natural to exclude the cases when $\tau\neq\inf\{s\ge0\dvtx
M_s\neq\hatM_s\}$. Note that for those cases the analogue of (\ref
{switch}) at $\inf\{s\ge0\dvtx  M_s\neq\hatM_s\}$ fails since this
is the first of the jump times after $\tau$, and the respective first
jump times of $M$ and $\hatM$ after $\tau$ will be different a.s.

%
\begin{prop}\label{propsplit2} Consider a pair of positive
nonincreasing pure jump processes $(M,\hatM)$, and suppose that $\tau
=\inf\{t\ge0\dvtx  M_t\neq\hatM_t\}\in(0,\infty)$ a.s. Then
$(M,\hatM)$ is a bifurcator with splitting time $\tau$ if and only if
the asymmetric Markov branching property of Definition~\ref{ambp}
holds relative to $\tau$, together with the binary junction conditions
(\ref{branch}) and (\ref{switch}).
\end{prop}

\begin{pf} The \emph{only if} part is straightforward. For the \emph{if}
part suppose that $(M,\hatM)$ satisfies the asymmetric Markov
branching property. Then Definition~\ref{ambp} implies that
$(M_t,\hatM_t,1_{\{\tau>t\}})$ is a Markov process in its natural
filtration. Furthermore, each component is Markovian in its own right
with respect to this filtration. More\vspace*{1pt} specifically, we see
that $M$ and $\hatM$ are fragmenters with some L\'evy measures
$\Lambda$ and $\widehat{\Lambda}$, and that $\tau$ is exponentially
distributed with some rate $\phi$. From (\ref{switch}), we have that
$\tau$ is a common jump time of all three components, so we can
consider $\tau$ as a marking time for each of the Poisson point
processes $((s,F_s),s\ge0)$ and $((s,\widehat{F}_s),s\ge0)$, with
mark 1 at $\tau$, say. By Lemma~\ref{454654156489784} below, this yields
marking kernels $K$ and $\widehat{K}$ from $(0,1)$ to $\{0,1\}$, from
which we define $p(u)=K(u,\{1\})$ and $\widehat{p}(u)=\widehat{K}(u,\{
1\})$. By standard results for marking and thinning Poisson point
processes, we find
\[
\mathbb{P}(F_\tau\in du)=\phi^{-1}p(u)\Lambda(du)\quad\mbox{and}\quad\mathbb{P}(\widehat{F}_\tau\in du)=\phi^{-1}
\widehat{p}(u)\widehat{\Lambda}(du),
\]
and the points before $\tau$, which by (\ref{branch}) are common to both
processes, have equal thinned intensity measures
\[
\bigl(1-p(u)\bigr)\Lambda(du)=\bigl(1-\widehat{p}(u)\bigr)\widehat
{\Lambda}(du).
\]
Together with (\ref{switch}), these equations are equivalent to (\ref
{hatlambda}), and we easily deduce that $(M,\hatM)$ is indeed a
bifurcator with splitting time $\tau$ in the sense of Definition~\ref{bifdef}.
\end{pf}

%
\begin{lm}\label{454654156489784} Consider a filtration $\mathcal{F}$, an
$\mathcal{F}$-Poisson point process $(F_t,t\ge0)$ with intensity
measure $\Lambda$ on $(0,1)$ and cemetery $1$, an $\mathcal
{F}$-stopping time $\tau$  such that $F_\tau\neq1$ a.s. and such that
conditionally given $\tau>t$, we have $(\tau-t,(F_{t+s},s\ge
0))\stackrel{d}{=}
(\tau,(F_s,s\ge0))$, for all $t\ge0$. Then there exists a
marking kernel $K$ from $(0,1)$ to $\{0,1\}$ such that for a Poisson
point process
$((\widetilde{F}_t,\widetilde{m}_t),t\ge0)$ with intensity measure
$\Lambda^+(du,dm)=K(u,dm)\Lambda(du)$ and for
$\widetilde{\tau}=\inf\{t\ge0\dvtx \widetilde{m}_t=1\}$, we have
$((F_t,t\ge0),\tau)\stackrel{d}{=}((\widetilde{F}_t,t\ge
0),\widetilde{\tau})$.
\end{lm}

We prove Lemma~\ref{454654156489784} in the \hyperref[appA]{Appendix}.

%
\begin{rem}\label{rem1} While a splitting density $f$ and a
switching probability function $p$ together uniquely\vspace*{1pt} identify the distribution
of a bifurcator, for two given splitting densities $f$ and $\widehat
{f}$, there may not be an associated bifurcator $(M,\hatM)$.
Looking ahead at Theorem~\ref{embed}, this will, in fact, be the
typical case. On the\vspace*{1pt} other hand, for fragmenters $M$ and
$\hatM$ that can be coupled to form a bifurcator, there are typically
many other couplings as different bifurcators. This can
be seen from (\ref{hatf}), which for each $u$ leaves some choice of
$p(u)$ and $p(1-u)$. We will see in Remark~\ref{rem2} that for
any choice with both $p(u)$ and $p(1-u)$ in $[0,1]$, equation (\ref
{hatf}) for $1-u$ instead of $u$, which appears to give a second equation
relating $p(u)$ and $p(1-u)$, will automatically hold if (and only if)
the two fragmenters can be embedded in the same
fragmentation process.
\end{rem}

%
\begin{ex}\label{ex7} An extreme example of a switching probability
function is $p(u)=1$ for $u<1/2$ and $p(u)=0$ for $u>1/2$. In\vspace*{1pt}
words: switch if the other block is bigger. We then obtain from (\ref
{hatf}) that $u\widehat{f}(u)=uf(u)+(1-u)f(1-u)$ for $u>1/2$
and $u\widehat{f}(u)=0$ for $u<1/2$. Note that in the context of
Remark~\ref{rem1}, there is only one bifurcator $(M,\hatM)$
which has a given $f$ and this associated $\widehat{f}$ as splitting densities.
\end{ex}

\subsection{Size-biased branching}\label{sizebiased}

The instance of the bifurcator construction of the previous section with
\[
p(u) = 1 - u,\qquad0<u<1, %
\]
is of special interest.
We then say that the bifurcator $(M,\hatM)$ is
\emph{derived from $M$ by size-biased branching}, and use the notation
$(M,M^*)$ instead of $(M,\hatM)$ to indicate this special construction.
Note that the ``size'' involved in the size-biasing is the size of the residual
factor $1-u$ associated with decrements of $M$ by a factor of $u$, that
is, the
size relative to the current value of $M$ of the fragment that splits.

We note the following corollaries to the results obtained in the
previous section.

%
\begin{cor}
\label{corstar}
If $M$ is a fragmenter with splitting density $f$ and Laplace exponent
\[
\Phi(\rho) = \int_0^1 \bigl(1 -
u^{\rho} \bigr) u f(u) \,du, %
\]
then $M^*$ derived from
$M$ by size-biased branching is a fragmenter with splitting density
%
\begin{equation}
\label{symmf} f^*( u ) = u f(u) + (1-u) f(1-u), \qquad0 < u < 1,
\end{equation}
and Laplace exponent
%
\begin{eqnarray}\label{phist}
\Phi^*(\rho) &=& \int_0^1 \bigl(1
- u^{\rho} \bigr) u f^*(u) \,du
\nonumber\\[-8pt]\\[-8pt]
&=& \Phi(\rho+ 1) - \Phi(\rho+1,\rho+1),\qquad\rho> 0,\nonumber
\end{eqnarray}
where $\Phi(\rho+1,\rho+1)$ is given by 
%
\[
\Phi(\rho+1,\rho+1) = \int_{0}^1 (1 -
u)^{\rho+ 1} u f(u) \,du. %
\]
Moreover, $M$ is then derived from $M^*$ by switching according to $p^*$,
where
%
\begin{equation}
\label{pstar} p^*(u) = \frac{ (1-u) f(1-u)} { f^* (u) }, \qquad0 < u < 1.
\end{equation}
\end{cor}

%
\begin{cor}\label{cor6}
In the setting of the previous corollary, the following conditions are
equivalent:
\begin{longlist}[(iii)]
\item[(i)] $f$ is symmetric: $f(u) = f(1-u)$ for all
$0<u<1$;
\item[(ii)] $f = f^*$;
\item[(iii)] $\Phi=\Phi^*$;
\item[(iv)] $p^*(u) = 1-u$ for all $0<u<1$;
\item[(v)] $M \stackrel{d}{=}M^*$;
\item[(vi)] $(M,M^*) \stackrel{d}{=}(M^*,M)$.
\end{longlist}
\end{cor}

Observe from (\ref{symmf}) that whatever the splitting density $f$ of
$M$, the splitting
density $f^*$ of $M^*$ is symmetric.
Thus the operation of passing from the law of $M$ to the law of $M^*$
by size-biased
branching is a kind of symmetrisation of laws of fragmenters
corresponding to the elementary symmetrisation of density functions defined
by formula (\ref{symmf}).
The operation is idempotent: $M^{**} \stackrel{d}{=}M^*$.
So we make the following definition:

%
\begin{defn}[(Symmetrised fragmenter)]
For a fragmenter $M$ with splitting density $f$
call the fragmenter $M^*$ with splitting density $f^*$ as
in (\ref{symmf}) the \emph{symmetrisation} of $M$.
\end{defn}

This notion of size-biased branching and symmetrisation can clearly
be extended to fragmenters whose splitting measures do not have a density,
as is achieved by formula (\ref{phist}) for Laplace exponents.
Note, however, that the probabilistic meaning in terms of size-biased
branching, and even the analytic fact that $\Phi^{**} = \Phi^{*}$ is very
much obscured by the Laplace exponents. Also, if $\Phi^*(\rho)=\int
_{(0,1)}(1-u^\rho)\Lambda^*(du)$, we obtain from
(\ref{hatlambda})
%
\begin{equation}
\label{pstargen}\Lambda^*(du)=u\Lambda(du)+u\overline{\Lambda
}(du)\quad\mbox{and}\quad p^*=d\overline{\Lambda}/d(\Lambda+\overline{\Lambda}).
\end{equation}
As can be seen in the symmetry discussion leading up to (\ref{ssphi})
below, this operation of symmetrisation of $M$ projects
the collection of laws of all fragmenters $M$ onto the collection of laws
of fragmenters $M^*$ which are canonically associated with a
binary homogeneous fragmentation process via the mass of a uniformly randomly
tagged fragment. This raises the question: exactly how
is a fragmenter $M$ with splitting density $f$ related to the
binary homogeneous fragmentation process with splitting density $f^*$?
We answer this question in Section~\ref{fragemb}, after development of
the necessary framework in Section~\ref{stexhom}.

\subsection{Exchangeable fragmentation processes and self-similar CRTs}\label{stexhom}

In our context, we can express Bertoin's \cite{Ber-hom,Berss}
definitions of homogeneous and self-similar fragmentations as follows.
For $\alpha\in\mathbb{R}$, we say that a family $\Pi_\alpha=(\Pi
_\alpha(t),t\ge0)$ of refining partitions in $\mathcal{P}=\mathcal
{P}_{\mathbb{N}}$ is an \emph{exchangeable $\alpha$-self-similar
fragmentation process} if both:
\begin{itemize}
\item$\Pi_\alpha$ is exchangeable in that the distribution of $\Pi
_\alpha$ is invariant under permutations of $\mathbb{N}$;
\item$\Pi_\alpha$ is a right-continuous strong Markov process whose
transition kernel satisfies the branching property: for all $t\ge0$,
$s\ge0$, conditionally given $\Pi_\alpha(t)=\{B_i,i\ge1\}$, the
partition $\Pi_\alpha(t+s)$ has the same distribution as the
partition of $\mathbb{N}$ with blocks
$B_i\cap\Pi_\alpha^{(i)}(|B_i|^{-\alpha} s)$, $i\ge1$, where the
$\Pi_{\alpha}^{(i)}$, $i\ge1$, are independent copies of $\Pi
_\alpha$.
\end{itemize}
Usually we consider $\Pi_\alpha(0)=\{\mathbb{N}\}$, and we exclude
the trivial case $\Pi_\alpha(t)=\{\mathbb{N}\}$ for all $t\ge0$.
Then $\Pi_\alpha(\infty):=\lim_{t\rightarrow\infty}\Pi_\alpha
(t)=\{\{1\},\{2\},\ldots\}$.
In the case $\alpha=0$, the linear time-changes of $\Pi_\alpha
^{(i)}$ by asymptotic frequencies $|B_i|^\alpha$ disappear; this case
is called
an \emph{exchangeable homogeneous fragmentation process}. Bertoin
\cite{Ber-hom} showed that the distribution of $\Pi=\Pi_0$ can be
expressed in terms of an exchangeable \mbox{$\sigma$-}finite intensity
measure $\kappa(d\Gamma)$ on $\mathcal{P}\setminus\{\{\mathbb{N}\}
\}$ via a L\'evy--It\^o-type decomposition into elementary splits of
blocks $B$ by $\Gamma=\{\Gamma_i,i\ge1\}$ into $B\cap\Gamma_i$,
$i\ge1$. The measure $\kappa$ admits an integral
representation
\[
\kappa(d\Gamma)=c\sum_{n\ge1}\delta_{\{\mathbb{N}\setminus\{i\},\{i\}\}}(d
\Gamma)+\int_{\mathcal{S}^\downarrow}\kappa_\mathbf{s} (d\Gamma)
\nu_{\mathrm{ranked}}(d\mathbf{s}),
\]
for an \emph{erosion coefficient} $c\ge0$ and a \emph{ranked
dislocation measure} $\nu_{\mathrm{ranked}}$ on $\mathcal{S}^\downarrow
\setminus\{(1,0,\ldots)\}$ satisfying $\int_{\mathcal{S}^\downarrow
}(1-s_1)\nu_{\mathrm{ranked}}(d\mathbf{s})<\infty$, and where $\kappa
_\mathbf{s}$ is
Kingman's paintbox governing exchangeable partitions with asymptotic
frequencies $\mathbf{s}\in\mathcal{S}^\downarrow:=\{(s_i)_{i\ge
1}\dvtx
s_1\ge s_2\ge\cdots\ge0,\sum_{i\ge1}s_i\le1\}$.
In
the \emph{binary} case, $\kappa(\Gamma\in\mathcal{P}\setminus\{
\mathbb{N}\}\dvtx \Gamma_1\cup\Gamma_2\neq\mathbb{N})=0$, this
representation can be written as
\[
\kappa(d\Gamma)=c\sum_{n\ge1}\delta_{\{\mathbb{N}\setminus\{i\},\{i\}\}}(d
\Gamma)+\frac{1}{2}\int_{(0,1)}\kappa_{(s,1-s)}(d
\Gamma)\nu(ds),
\]
for a \emph{symmetric dislocation measure} $\nu$ on $(0,1)$
satisfying $\int_{(0,1)} s(1-s)\nu(ds)<\infty$ and $\nu=\overline
{\nu}$, where $\overline{\nu}$ is the push-forward of $\nu$ under
$u\mapsto1-u$, so
$\nu(du)=\nu_{\mathrm{ranked}}(s_1\in du)+\nu_{\mathrm{ranked}}(s_2\in du)$
and $\nu_{\mathrm{ranked}}(s_1\in\cdot)=\nu( \cdot\cap(\frac
{1}{2},1))+ \frac{1}{2}\nu( \cdot\cap\{\frac{1}{2}\})=\overline
{\nu}( \cdot\cap(\frac{1}{2},1))+\frac{1}{2}\overline{\nu}(
\cdot\cap\{\frac{1}{2}\})$.

We denote by $|A_\alpha^n(t)|$ the asymptotic frequency of the block
$A_\alpha^n(t)$ of $\Pi_\alpha(t)$ containing $n$. For $\alpha=0$,
the process $|A_0^n(t)|$ is a fragmenter, and $\xi_n(t)=-\log
|A_0^n(t)|$ has Laplace exponent
%
\begin{equation}
\label{ssphi}\Phi^*(\rho)=c+c\rho+\int_{(0,1)}
\bigl(1-u^\rho\bigr)u\nu(du);
\end{equation}
see \cite{Ber-hom}. Self-similar and homogeneous fragmentation
processes are pathwise related by nonlinear time-change \cite{Berss}.
Specifically,
%
\begin{eqnarray}\label{afreq}
\bigl|A_\alpha^n(t)\bigr|=\exp\bigl(-\xi_n\bigl(
\eta_n(t)\bigr) \bigr)
\nonumber\\[-8pt]\\[-8pt]
\eqntext{\displaystyle\mbox{where }\eta_n(t)=\inf
\biggl\{u\ge0\dvtx \int_0^ue^{-\alpha\xi_n(w)}\,dw>t
\biggr\},}
\end{eqnarray}
is a self-similar Markov process, and for all $\alpha>0$, we have
$|A_\alpha^n(t)|=0$ for $t\ge\int_0^\infty e^{-\alpha\xi_n(w)}\,dw$.

It was shown by Haas and Miermont \cite{HM04} that for every
exchangeable self-similar
fragmentation process $\Pi_\alpha$ with index $\alpha> 0$, zero
erosion $c=0$ and infinite dislocation measure $\nu_{\mathrm{ranked}}$ on
$\mathcal{S}_1^\downarrow=\{\mathbf{s}\in\mathcal{S}^\downarrow
\dvtx
\sum_{i\ge1}s_i=1\}$, there is an
\mbox{associated} compact continuum random tree $(\mathcal{T},d,0,\mu)$.
Vice versa, such a continuum random tree (CRT) allows an embedding of a
self-similar fragmentation process. Specifically, a CRT is a random
weighted and rooted $\mathbb{R}$-tree. A weighted and rooted $\mathbb
{R}$-tree $(T,d,0,\mu)$ is a complete, separable metric space $(T,d)$
with a \emph{root} $0\in T$ and a probability measure $\mu$ on the
Borel sets of $(T,d)$, such that the following tree property holds:
\begin{itemize}
\item Any two points $\sigma,\sigma^\prime\in T$ are connected by a
unique injective path $[\![\sigma,\sigma^\prime]\!]$. Furthermore, this
path can be uniquely parametrised $[\![\sigma,\sigma^\prime]\!]=\{
g_{\sigma,\sigma^\prime}(t),0\le t\le d(\sigma,\sigma^\prime)\}$
by an isometry $g_{\sigma,\sigma^\prime}\dvtx [0,d(\sigma,\sigma
)]\rightarrow T$ with $g_{\sigma,\sigma^\prime}(0)=\sigma$ and
$g_{\sigma,\sigma^\prime}(d(\sigma,\sigma^\prime))=\sigma^\prime$.
\end{itemize}
We also write $]\!]\sigma,\sigma^\prime]\!]:=[\![\sigma,\sigma^\prime
]\!]\setminus\{\sigma\}=\{g_{\sigma,\sigma^\prime}(t),0<t\le
d(\sigma,\sigma^\prime)\}$.

When there is no ambiguity about $d$, $0$ or even $\mu$, we simply
write $(\mathcal{T},\mu)$ or even $\mathcal{T}$ to refer to a CRT
$(\mathcal{T},d,0,\mu)$. For the purpose of convergence of compact
weighted and rooted $\mathbb{R}$-trees, we will identify $(T,d,0,\mu
)$ and $(T^\prime,d^\prime,0^\prime,\mu^\prime)$ if there is an
isometry between $(T,d)$ and $(T^\prime,d^\prime)$ that maps $0$ to
$0^\prime$ and pushes $\mu$ forward to $\mu^\prime$. The set
$\mathbb{T}$ of such isometry classes can then be equipped with the
so-called Gromov--Hausdorff--Prohorov distance $d_{\mathrm{GHP}}$. Then
$(\mathbb{T},d_{\mathrm{GHP}})$ is a Polish metric space. See, for example,
\cite{Eva-10} for further background on the space $(\mathbb{T},d_{\mathrm{GHP}})$.

For an $\mathbb{R}$-tree $(T,d,0,\mu)$, let $T^t:=\{\sigma\in
T\dvtx  d(0,\sigma)>t\}$, $t\ge0$, and define fringe subtrees
$T_\sigma:= \{\sigma^\prime\in T\dvtx \sigma\in
[\![0,\sigma^\prime]\!]\}$, $\sigma\in T$. An $\alpha
$-self-similar CRT is a random weighted and rooted $\mathbb{R}$-tree
$(\mathcal{T},d,0,\mu)$, or its isometry class with distribution on
the Borel space of $(\mathbb{T},d_{\mathrm{GHP}})$, such that:
\begin{itemize}
\item$\mu$ is nonatomic with dense support a.s., $\mu(\mathcal
{T}_\sigma)>0$ for all $\sigma\in\mathcal{T}$ with $\mathcal
{T}_\sigma\neq\{\sigma\}$, while $\mu([\![0,\sigma]\!])=0$ for all
$\sigma\in\mathcal{T}$, and
\item for all $t\ge0$, the connected components $(\mathcal
{T}_i^t,i\ge1)$ of $\mathcal{T}^t$, completed by a root vertex $0_i$,
are such that given $(\mu(\mathcal{T}_i^t),i\ge1)=(m_i,i\ge1)$, for
some $m_1\ge m_2\ge\cdots\ge0$, the trees
\[
\bigl(\mathcal{T}_i^t,m_i^{-\alpha}\,d\bigl|_{\mathcal{T}_i^t},0_i,m_i^{-1}
\mu\bigr|_{\mathcal{T}_i^t}\bigr),\qquad i\ge1
\]
are independent and identically distributed isometric copies of
$(\mathcal{T},d,0,\mu)$.
\end{itemize}
Recently, Stephenson \cite{Ste-13} extended this class by relaxing the
first bullet point to
allow a support that is not dense, atoms of $\mu$, and/or positive
weights on branches, so as to include all dislocation measures $\nu
_{\mathrm{ranked}}$ and erosion $c>0$. The CRT $(\mathcal{T},\mu)$
constructed in \cite{HM04,Ste-13} is such that
%
\begin{equation}
\label{embedpart}\qquad \Pi_\alpha^*(t)=\bigl\{\bigl\{j\ge1\dvtx
\Sigma_j^*\in\mathcal{T}^t_i\bigr\},i\ge1
\bigr\}\cup\bigl\{\{j\},j\ge1\dvtx \Sigma_j^*\notin
\mathcal{T}^t\bigr\},\qquad t\ge0
\end{equation}
has the same distribution as $(\Pi_\alpha(t),t\ge0)$, where
conditionally given $(\mathcal{T},\mu)$ the sequence $\Sigma_n^*$,
$n\ge1$,
is independent and identically distributed according to~$\mu$. It was
shown in \cite{HM04,Ste-13} that the subtrees $\mathcal
{R}_k^*=\bigcup_{j=1}^k[\![0,\Sigma_j^*]\!]\subset\mathcal{T}$ converge
a.s. in the Hausdorff sense for embeddings in $\ell_1(\mathbb{N})$,
and this easily entails $d_{\mathrm{GHP}}((\mathcal{R}_k^*,\mu
_k^*),(\mathcal{T},\mu))\rightarrow0$ a.s., as $k\rightarrow\infty
$, where $\mu_k^*$ is the push-forward of~$\mu$ under the projection
map $\pi^{\mathcal{R}_k^*}\dvtx \mathcal{T}\rightarrow\mathcal
{R}_k^*$, $u\mapsto g_{0,\sigma}(\sup\{t\ge0\dvtx  g_{0,\sigma
}(t)\in\mathcal{R}_k^*\})$.
Also, $\mu$ is then recovered in accordance with Aldous's theory of
consistent leaf-exchangeable families of trees $(\mathcal{R}_k^*,k\ge
1)$ as the weak limit of the uniform distribution $\nu_k^*$ on $\Sigma
_1^*,\ldots,\Sigma_k^*$, as $k \rightarrow\infty$.

\subsection{Embedding fragmenters in homogeneous binary fragmentation processes}\label{fragemb}

Let $\Pi= (\Pi(t), t \ge0 )$ be a binary homogeneous fragmentation
process starting from $\Pi(0)=\{\mathbb{N}\}$, with absolutely
continuous symmetric dislocation measure $\nu(du)=f^*(u)\,du$, for some
\emph{symmetric splitting density} $f^*$ on $(0,1)$,
so $\Pi$ takes values in the set of partitions of $\mathbb{N}$.
Let $A = (A_t, t \ge0 )$ be a process with values in subsets of
$\mathbb{N}$.
Call $A$ a \emph{path in $\Pi$} if:
\begin{itemize}
\item$A_t \in\Pi(t)$ for all $t \ge0$;
\item$A_t$ is decreasing in the inclusion partial order, as $t$ increases.
\end{itemize}
%

\begin{defn}[(Markovian path)]
We call a path $A$ a \emph{Markovian
path in~$\Pi$}~if:
\begin{itemize}
\item$(A_t, t \ge0 )$ is adapted to some filtration $(\mathcal{F}_t,
t \ge0 )$
with respect to which $(\Pi(t), t \ge0)$ is Markovian, in such a way that
for each $s \ge0$ the process $((A_t,\Pi(t)\cap A_s), t \ge s )$ and
the restriction of $(\Pi(t), t \ge s)$ to $\mathbb{N}- A_s$ are conditionally
independent given $\mathcal{F}_s$.
\end{itemize}
\end{defn}
%

To explain the terminology, think of $\Pi$ as embedded by suitable
time change in the $\alpha$-self-similar
continuum random tree (CRT) associated with $\Pi$ for $\alpha> 0$ by
Haas and Miermont \cite{HM04}, or by Stephenson \cite{Ste-13} when
$\nu$ is finite.
Then $A_t$ represents the set of leaf labels above some internal vertex
$v_t$ of the CRT in some path leading from the root
to a leaf vertex $v_\infty$ of the tree.

Let $M$ be a fragmenter. Say that $M$ admits a \emph{Markovian embedding
in $\Pi$}
if it is possible to construct $\Pi$ jointly with a Markovian path $A$
such that
\[
(M_t, t \ge0 ) \stackrel{d} {=}\bigl(|A_t|, t \ge0\bigr),
\]
where $|A_t|$ is the asymptotic frequency of $A_t$, which is known to
exist almost surely, simultaneously for all $t \ge0$ and all
sets $A_t \in\Pi(t)$.
In terms of an associated CRT construction, the jumps of
$|A|:=(|A_t|,t\ge0)$ would then describe the spinal partition of mass
in the CRT along a
spine leading from the root to some random leaf of the CRT.
The most basic example is provided by the next proposition:

%
\begin{prop}\label{symmbif}
For each positive integer $n$, let $A^{n}_t$ be the block of $\Pi(t)$
containing $n$. Then:
\begin{longlist}[(ii)]
\item[(i)] (\emph{Bertoin} \cite{Ber-hom})
$A^n=(A^n_t,t\ge0)$ is a Markovian path in $\Pi$ such that $|A^n|$ is
a fragmenter with splitting density $f^*$.
\item[(ii)] For each pair of positive integers $m$ and $n$,
the pair $(|A^n|,|A^m|)$ is a symmetric bifurcator,
each derived by size-biased branching from the other.
\end{longlist}
\end{prop}

\begin{pf}
Part (ii) follows easily from the (strong) homogeneous branching property
of $\Pi$ and Corollary~\ref{corstar}.
\end{pf}

Due to the natural embeddings provided by this proposition, we call the
fragmenter $M^*$ with splitting density $f^*$ the
\emph{canonical fragmenter} associated with $\Pi$.
In terms of the corresponding mass fragmentation $(|\Pi(t)|^\downarrow,t\ge0)$ of asymptotic frequencies $|\Pi_i(t)|$, $i\ge1$, ranked in
decreasing order $|\Pi(t)|^\downarrow\in\mathcal{S}^\downarrow$,
the process $M^*$ describes the
evolution of the mass of the fragment containing a randomly tagged
particle. We state more generally without assuming the existence of the
densities $f$ and $f^*$:

%
\begin{teo}
\label{embed}
Every fragmenter $M$ with Laplace exponent $\Phi(\rho)=\break\int
_{(0,1)}(1-u^\rho)\Lambda(du)$ admits a Markovian embedding in an
exchangeable binary homogeneous fragmentation process $\Pi$. The
distribution of $\Pi$ is unique. Its symmetric dislocation measure is
given by $\nu=\Lambda+\overline{\Lambda}$. The canonical
fragmenter of $\Pi$ is the symmetrisation of $M$.
\end{teo}

This confirms points 1--3 of the \hyperref[sec1]{Introduction}. In the absolutely
continuous case, we can rephrase as follows. $M$ can be embedded in
$\Pi$ if and only if the splitting density of $\Pi$ is $f^*$, the
symmetrisation of the splitting density $f$ of $M$.

%
\begin{rem}\label{rem2}
Consider\vspace*{1pt} any bifurcator $(M,\hatM)$. By
adding equation (\ref{hatf}) and the equation we obtain by
substituting $u$ by $1-u$ in (\ref{hatf}), we see that $f$ and
$\widehat
{f}$ have the same symmetrisation $f^*$. By Theorem
\ref{embed}, $M$ and $\hatM$ can each be embedded in the same binary
homogeneous fragmentation process $\Pi$. In fact, the argument used to
prove this theorem, can be adapted to prove that the bifurcator admits
an embedding in $\Pi$. We leave the details of this to
the reader.
\end{rem}

To prepare for the proof of the theorem, we start with some remarks about
paths $A$ in $\Pi$.
For $t \ge0$, let $N_t:= \min A_t$. Clearly, $N_0=1$.
The fact that $A_t$ decreases as $t$ increases implies that $(N_t,t \ge
0)$ is some increasing process.
Furthermore, $A_t\in\Pi(t)$ implies that
%
\begin{equation}
\label{switching} A_t = A^{N_t}_t, \qquad t \ge0.
\end{equation}
Assuming for simplicity that $N_t$ tends to $\infty$ as $t \to\infty$,
let $0= \tau_0 < \tau_1 < \tau_2 < \cdots$ be the successive times
of jumps of $(N_t, t \ge0)$,
and set $N(n) = N_{\tau_n}$, $n\ge0$.
Then $A_t = A_t^{N(n)}$ for $t \in[\tau_{n}, \tau_{n+1})$. Note that
given the random sequence $1 = N(0) < N(1) < N(2) < \cdots$,
the times $\tau_n$ can be recovered without further reference to~$A$,
from the family of paths $A^n$ associated with
$\Pi$, as $\tau_{n} = \inf\{t\ge0\dvtx  A^{N(n-1)}_t \ne A^{N(n)}_t
\}$ for each $n \ge1$.
Thus there is a natural correspondence between paths $A$ in $\Pi$ and
increasing random sequences $(N(n), n \ge1)$ subject
to the constraint that $N(n) = \min A_{\tau_n}^{N(n)}$, where in
general, the possibility of a finite increasing sequence of random
length must also be allowed.

In connection with the $\alpha$-self-similar CRT $\mathcal{T}$
derived from $\Pi$, notice that the random times $\tau_n$ are defined
in a way which allows corresponding random times
\[
\tau_{n,\alpha}:= \int_0^{\tau_n}
\bigl|A^{N(n-1)}_t\bigr|^{\alpha} \,dt
\]
to be defined, and that in $\mathcal{T}$ there is a junction
vertex $V_n$ at height $\tau_{n,\alpha}$ at which the paths to leaves
labelled $N(n-1)$ and $N(n)$ diverge. Here it is assumed that
the CRT is equipped with a random sample $\Sigma_1^*, \Sigma_2^*,
\ldots$ of its leaves according to its mass measure, and that the
homogeneous fragmentation
$\Pi= \Pi_{0}$ has been constructed by time changing the $\alpha
$-self-similar fragmentation associated with the tree and the leaf sample.
The path from the root $0$ to $V_n$ in the CRT then passes through
$V_1, \ldots, V_{n-1}$. Since the CRT is compact, we
find a convergent subsequence of $(V_n,n\ge1)$ with limit $\Sigma$,
say. Because $(V_n,n\ge1)$ is increasing for the genealogical partial
order $\prec$ that puts $\sigma\prec\sigma^\prime$ if and only if
$\sigma\in[\![0,\sigma^\prime[\![$, the sequence converges to the same
limit. Note that $\Sigma$ must be a leaf almost surely
because if $\Sigma$ is not a leaf, then the fringe subtree $\mathcal
{T}_\Sigma$ at $\Sigma$ will have positive mass, but then $A_\infty
=\bigcap_{t\ge0}A_t$ will have positive limiting frequency, which
contradicts $\Pi(\infty)=\{\{1\},\{2\},\ldots\}$.
The path from $0$ to $\Sigma$ in $\mathcal{T}$ starts by following
the path to~$1$, then branches off in the
direction of $\Sigma_{N(1)}$, then branches again in the direction of
$\Sigma_{N(2)}$, etc. This could be formalised to give a one-to-one
correspondence between paths in $\Pi$ and paths in $\mathcal{T}$.


\begin{pf*}{Proof of the existence part of Theorem~\ref{embed}}
For a fragmenter with Laplace exponent $\Phi(\rho)=\int_{(0,1)}(1-u^\rho
)\Lambda(du)$, consider its symmetrisation $\Phi
^*(\rho)=\int_{(0,1)}(1-u^\rho)\Lambda^*(du)$ with $\Lambda^*$
given in (\ref{pstargen}). By (\ref{ssphi}), this
is the Laplace exponent of the canonical fragmenter of a fragmentation
process $\Pi$ with symmetric dislocation measure $\nu
(du)=u^{-1}\Lambda^*(du)=\Lambda(du)+\overline{\Lambda}(du)$.
Consider the canonical fragmenter $M^1=|A^1|$ obtained from the blocks
$A_t^1$, $t\ge0$, of $\Pi$ containing $1$ as in Proposition~\ref
{symmbif}. Construct a bifurcator $(M^1,M)$ by
switching from $M^1$ according to $p^*$ in (\ref{pstargen}), as needed to
create $M$ up to some branching time $\tau_1$.
For $0 \le t < \tau_1$ let $A_t = A_t^1$, and let
$A_{\tau_1} = A_{\tau_1-}^1\setminus A_{\tau_1}^1$, that is, the
block that splits off
from the block containing $1$ at time $\tau_1$.
To continue the construction of $A$ after time $\tau_1$, apply the
strong Markov property of $\Pi$ at $\tau_1$,
and let $N(1) = \min(A_{\tau_1})$. Set $A_t = A_t^{N(1)}$ for $ \tau
_1 \le t < \tau_2$ where
$\tau_2 - \tau_1$ is the branching time between
$M^{(1)}:= (|A_{\tau_1 + s }^{N(1)}|/|A_{\tau_1}|, s \ge0)$, which
is another
copy of $M^1$, and some further copy of $M$ created conditionally
given $M^1$ and $M^{(1)}$ by the same rule.
It is clear that continuing like this creates time segments $\tau_i -
\tau_{i-1}$
which are independent and identically distributed, and fresh copies of
$M^1$ as needed. The process $A$ with the desired feature that $|A|
\stackrel{d}{=}M$
can be created for all times $t \ge0$. Moreover, by construction this process
$A$ is a Markovian path in $\Pi$.
\end{pf*}

We postpone the uniqueness part of the proof. Specifically, points 1
and 3 of the \hyperref[sec1]{Introduction} have now been proved, while point 2 is
postponed to Section~\ref{genline}, where we first establish points
4--6 in Theorem~\ref{teogrowth}.

\subsection{Mass distributions}\label{sectmass}

The bifurcator $(M,M^*)$ with $M^*$ derived from $M$ by size-biased branching
plays a key role in following discussions.
This section collects together
some basic formulae for the joint distribution of the branching time
$\tau: = \inf\{t\ge0\dvtx  M_t \ne M_t^* \}$ and the
decrements
\[
M_{\tau-} - M_\tau= M_\tau^*\quad\mbox{and}\quad
M_\tau= M_{\tau-}^* - M_\tau^*,
\]
where $M_{\tau-} = M^*_{\tau-}$.
The triple of nonnegative variables
$(1 - M_{\tau-}, M_{\tau^*}, M_\tau)$ with sum $1$ is of special interest.
In a suitably defined random $\mathbb{R}$-tree $\mathcal{R}_{\Sigma,\Sigma^*}$ spanned by a root $0$ and two
leaves $\Sigma$ and $\Sigma^*$,
this triple represents the masses of three connected components of the tree
formed by removal of a particular random junction vertex of the tree.
As indicated in the previous section, this subtree $\mathcal
{R}_{\Sigma,\Sigma^*}$ may be naturally embedded in
a self-similar CRT $\mathcal{T}$ associated with a fragmentation
process, whose canonical fragmenter is $M^*$.
The joint distribution of this triple is determined by a formula for its
joint moments provided by Gnedin and Pitman \cite{gp}, page 477, where this
triple is denoted $(G,H,D)$, with the following more elementary interpretation:
$(G,1-D)$ is the interval component covering $U$ in the complement of the
range of $(M_t, t \ge0 )$, for $U$ a uniform$(0,1)$ variable independent
of $M$, and $H = 1 - D-G$ is the length of this interval component.
The following formulae can be read either from the discussion of the
previous sections, or from
\cite{gp}, page 477.

Recall first that for $M$ with L\'evy exponent $\Phi$, L\'evy measure
$\lambda(x)\,dx$,
and splitting density $f$, the branching time $\tau$ has exponential
distribution with
rate
\[
\Phi(1) = \int_0^1 u ( 1 - u ) f(u) \,du = \int
_0^\infty\bigl( 1 - e^{-x} \bigr)
\lambda(x)\,dx. %
\]
The process $(M_t, 0 \le t < \tau)$ is then the negative exponential
of a killed subordinator with L\'evy measure $e^{-x} \lambda(x)\,dx$
and killing at rate $\Phi(1)$.

Thus
\[
\mathbb{E}\bigl[ M_t^\rho\giv\tau> t \bigr] = \mathbb{E}
\bigl[ M_{\tau-}^{\rho} \giv\tau= t \bigr] = e^{- t \Phi_0(\rho) },
\]
where
\[
\Phi_0(\rho) = \int_0^1 \bigl(1 -
u^\rho\bigr) u^2 f(u) \,du = \Phi(\rho+ 1 ) - \Phi(1)
\]
and hence by conditioning on $\tau$
%
\begin{equation}
\mathbb{E}\bigl[ M_{\tau-}^\rho\bigr] = \frac{\Phi(1) }{\Phi(1) + \Phi
_0(\rho) } =
\frac{\Phi(1) }{\Phi(\rho+ 1)}.
\end{equation}
This is \cite{gp}, formula (57) or (59).
From \cite{gp}, formula (28), or from (\ref{vdens}),
%
\begin{equation}
\mathbb{P} \biggl( \frac{M_\tau}{M_{\tau-}} \in du \biggr) = \bigl(\Phi(1)
\bigr)^{-1}( 1 - u)u f(u) \,du,
\end{equation}
hence
%
\begin{equation}
\mathbb{E} \biggl[ \biggl( \frac{M_\tau}{M_{\tau-}} \biggr)^\rho\biggr
] =
\frac{ \Phi(\rho+ 1 ) - \Phi(\rho) } {\Phi(1) }.
\end{equation}
Moreover, $ M_{\tau-}$ and $M_\tau/M_{\tau-}$ are independent, so
the last two formulae combine to give
%
\begin{equation}
\mathbb{E}\bigl( M_\tau^\rho\bigr) = \frac{ \Phi(\rho+ 1 ) - \Phi(\rho) }{
\Phi(\rho+ 1) },
\end{equation}
which is a simplification of \cite{gp}, formula (60).
Next, $M_\tau^*:= M_{\tau-} - M_\tau$ is a
size-biased pick from the decrements of $M$, whence
%
\begin{eqnarray}
\mathbb{E}\bigl[ \bigl(M_\tau^*\bigr)^\rho\bigr] =
\frac{\Phi(\rho+ 1, \rho+
1)}{\Phi(\rho+ 1)}
\nonumber\\[-8pt]\\[-8pt]
\eqntext{\displaystyle\mbox{where } \Phi(\rho+ 1, \rho+ 1 ) = \int
_0^1( 1 - u )^{\rho+ 1 } u f(u) \,du.}
\end{eqnarray}
Note that for positive integers $\rho= n$ say, this
is a
linear combination of
evaluations of $\Phi(k)$ at integers $k \le n + 1$,
as indicated in \cite{gp}, formula (25).
In principle, these Mellin transforms determine the distributions of
$M_\tau$ and $M_\tau^*$, but there do not seem to be simple formulae for
the densities of these
variables except in special cases.

Observe that the expected masses of the three components in the junction
split are
%
\begin{eqnarray}
\mathbb{E}\bigl( M_\tau^* \bigr) &=& \bigl(2 \Phi(1) - \Phi(2)\bigr)/
\Phi(2),
\\
\mathbb{E}( 1 - M_ {\tau- } ) &=& \mathbb{E}( M_ {\tau} ) = \bigl(
\Phi(2) - \Phi(1) \bigr)/\Phi(2).
\end{eqnarray}
It does not seem obvious intuitively why the expectations
of $( 1 - M_{\tau- } )$ and $M_{\tau}$ are always equal.

Recall from (\ref{phist}) that $\Phi^*(\rho):= \Phi(\rho+ 1 ) -
\Phi
(\rho+1, \rho+ 1)$
is the Laplace exponent corresponding to $M^*$ with splitting density
$f^*(u) = uf(u) + (1-u) f(1-u)$. So we obtain the following extension
of Corollary~\ref{cor6}:

\begin{cor} Each of the following two conditions is also equivalent to
the symmetry of $f$, which we characterised in Corollary~\ref{cor6} as
$f = f^*$, and as $\Phi= \Phi^*$:%
\begin{longlist}[(viii)]
\item[(vii)] $M_{\tau} /M_{\tau-}\stackrel{d}{=}M_{\tau} ^*
/M_{\tau
-}$;
\item[(viii)] $M_{\tau}\stackrel{d}{=}M_{\tau}^*$.
\end{longlist}
\end{cor}
%

\subsection{Edge lengths and exponential functionals}\label{sectlength}

Continuing to suppose that $(M,M^*)$ is the bifurcator derived from
$M$ by size-biased branching, as well as the basic triple of masses
$(1 - M_{\tau-}, M_{\tau}^*, M_\tau)$ with sum $1$, for each $\rho>0$
we may consider the triple of \emph{exponential functionals}
\[
L_0 = \int_0^\tau
M_{t}^\rho \,dt, \qquad L_\Sigma= \int
_\tau^\infty M_{t}^\rho \,dt\quad\mbox{and}\quad L_* = \int_\tau^\infty
\bigl(M_{t}^*\bigr)^\rho \,dt, %
\]
which can be interpreted as the lengths of branches in a suitably defined
random \mbox{$\rho$-}self-similar $\mathbb{R}$-tree with three branches
meeting at a junction point,
these branches being labelled by $0 $ for the root, and $\Sigma$ and
$*$ for
the two leaves associated with $M$ and $M^*$, respectively.
Note that the definition of the $L_i$ depends on the parameter $\rho$,
which is suppressed in the notation.
In particular, suppose that the L\'evy measure satisfies the regular variation
condition $\int_x^\infty\lambda(x)\,dx=x^{-\alpha}\ell(1/x)$ as
$x\downarrow0$, for some $\alpha\in(0,1)$ and some function $\ell
\dvtx (0,\infty)\rightarrow(0,\infty)$ that is slowly varying at
$\infty$. Then the above functionals for
$\rho= \alpha$ are of special interest \cite{GPYalpha}. They govern
the asymptotics
of how numbers of new branch points grow along the three branches, as
new branch points are selected by size-biased sampling from the mass
distribution
on the $\mathbb{R}$-tree $\mathcal{R}_{\Sigma,\Sigma^*}$, which
assigns the decrements of $M$
except
$M_{\tau-}-M_\tau= M_\tau^*$
to the branch from the root $0$ to leaf $\Sigma^*$, and the
remaining decrements of $M^*$ to the branch from the junction point to
leaf $\Sigma^*$.
Let
\begin{eqnarray*}
L_{0 \Sigma}&:=& L_{0 } + L_\Sigma= \int
_0^\infty M_t^\rho \,dt = \int
_0^\infty e^{- \rho\xi_t } \,dt,
\\
L_{0 *} &:=& L_{0 } + L_* = \int_0^\infty
\bigl(M_t^*\bigr)^\rho \,dt = \int_0^\infty
e^{- \rho\xi_t ^*} \,dt,
\end{eqnarray*}
where $\xi$ and $\xi^*$ are the two subordinators associated with $M$
and $M^*$.
According to a known formula for subordinators \cite{BY05},
%
\begin{equation}
\label{prodmom} \mathbb{E}\bigl( L^n_{0 \Sigma} \bigr)
= \frac{n ! }{\Phi(\rho) \cdots\Phi(n \rho) } \qquad\mbox{for all $n\in
\mathbb{N}$,}
\end{equation}
where $\Phi$ is the Laplace exponent of $\xi$,
and the same holds for $L_{0 *}$ instead of
$L_{0 \Sigma}$ with the Laplace exponent $\Phi^*$ of $\xi^*$ instead of
$\Phi$.
Now
\[
L_0 = \int_0 ^\tau e^{- \rho\xi_t }
\,dt = \int_0^\infty e^{- \rho
\eta_t } \,d t,
\]
where $\eta_t = \xi_t 1_{\{\tau> t\}}+ \infty1_{\{\tau\le t\}}$ is
another subordinator, whose L\'evy measure is
$e^{-x} \lambda(x)\,dx + \Phi(1) \delta_\infty(dx)$ for $\lambda
(x)\,dx$ the
L\'evy measure of $\xi$. It follows that the Laplace exponent of $\eta
$ at $\rho$
is $\Phi_\eta(\rho)=\Phi(1)+\int_0^\infty(1-e^{-\rho
x})e^{-x}\lambda(x)\,dx=\Phi(1)+\int_0^\infty(1-e^{-(\rho
+1)x})\lambda(x)\,dx-\int_0^\infty(1-e^{-x})\lambda(x)\,dx=\Phi(\rho+ 1)$,
and hence that
%
\begin{equation}
\label{lomoms} \mathbb{E}\bigl( L_0 ^n \bigr) =
\frac{n!}{\Phi_\eta(\rho)\cdots\Phi_\eta
(n\rho)}= \frac{n ! }{\Phi(\rho+1 ) \cdots\Phi(n \rho+ 1) }\qquad\mbox
{for all $n\in\mathbb{N}$.}\hspace*{-20pt}
\end{equation}
Moments of $L_\Sigma$ and $L_*$ can now be found using the distributional
identities $L_\Sigma\stackrel{d}{=}M_\tau^\rho\widehat{L}_{0
\Sigma}$ and
$L_* \stackrel{d}{=}(M_\tau^*)^\rho\widehat{L}_{0 *}$ where
$\widehat{L}_{0 \Sigma}$ is independent of $M_\tau$ with $\widehat
{L}_{0 \Sigma} \stackrel{d}{=}L_{0 \Sigma}$,
and $\widehat{L}_{0 *}$ is independent of $M_\tau^*$ with
$\widehat{L}_{0 *} \stackrel{d}{=}L_{0 *}$.
Thus
%
\begin{eqnarray}\label{lmoremoms}
\quad \mathbb{E}\bigl( L_\Sigma^n \bigr) &=&
\mathbb{E}\bigl(M_\tau^{\rho n } \bigr) \mathbb{E}
\bigl(L_{0 \Sigma}^{n } \bigr) = \frac{\Phi(\rho n+1)-\Phi(\rho n)}{\Phi
(\rho n+1)}
\frac{n!}{\Phi(\rho)\cdots\Phi(n\rho)},
\\
\label{lmoremoms1} \mathbb{E}\bigl( L_*^n \bigr) &=& \mathbb{E}\bigl(
\bigl(M_\tau^*\bigr)^{\rho n } \bigr) \mathbb{E}
\bigl(L_{0 *}^{n } \bigr) = \frac{\Phi(\rho n+1,\rho n+1)}{\Phi(\rho
n+1)}
\frac{n!}{\Phi^*(\rho)\cdots\Phi^*(n\rho)}.
\end{eqnarray}
Note the two identities in distribution
%
\begin{equation}
\label{lids} L_{0 \Sigma} \stackrel{d} {=}L_0 +
M_\tau^{\rho} \widehat{L}_{0
\Sigma}\quad\mbox{and}\quad L_{0 *} \stackrel{d} {=}L_0 + \bigl(M_\tau^*
\bigr)^{\rho} \widehat{L}_{0 *},
\end{equation}
where
$L_{0 \Sigma} \stackrel{d}{=}\widehat{L}_{0 \Sigma}$
with
$\widehat{L}_{0 \Sigma}$ independent of $(L_0, M_\tau)$, and
$
L_{0 *}
\stackrel{d}{=}
\widehat{L}_{0 *}
$
with
$\widehat{L}_{0 *}$ independent of $(L_0, M_\tau^*)$.
As checks, the two equalities of means implied by (\ref{lids})
are easily seen to be consistent with previous formulae.
The equalities of higher moments in (\ref{lids}) provide
identities involving joint moments such as
$\mathbb{E}( L_0 ^j M_\tau^{k \rho} )$ for positive integers
$j$ and $k$. In particular, $\mathbb{E}( L_0 M_\tau^{\rho} )$ is
determined by the second moment formula. But the third moment
formula only gives access to a linear combination of
$\mathbb{E}( L_0 M_\tau^{2 \rho} )$ and $\mathbb{E}( L_0 ^2 M_\tau
^{\rho} )$,
which is not so useful.

\section{Bead splitting processes and continuum random trees}\label{genline}

%

Recall from Section~\ref{stexhom} that every self-similar CRT
$(\mathcal{T},\mu)$ gives rise to a growing family $(\mathcal
{R}_k^*,\mu_k^*)$ of
weighted $\mathbb{R}$-trees that converges to $(\mathcal{T},\mu)$.
As we will demonstrate more formally below, picking $\Sigma_{k+1}^*$
from $\mu$ means that a junction point $J_k^*\in\mathcal{R}_k^*$ is
picked from $\mu_k^*$ and that $\Sigma_{k+1}^*$ is then picked within a
subtree rooted at $J_k^*$, which is a rescaled copy of $(\mathcal
{T},\mu)$, by self-similarity. Then
\[
\mathcal{R}_{k+1}^*=\mathcal{R}_k^*\cup\,\bigl]\!\bigl]J_k^*,
\Sigma_{k+1}^*\bigr]\!\bigr]\quad\mbox{and}\quad\mu_{k+1}^*=
\mu_k^*-\mu_k^*\bigl(\bigl\{ J_k^*\bigr\}
\bigr)\delta_{J_k^*}+\mu_{k,k+1}^*,
\]
where $\mu^*_{k,k+1}$ denotes the projection onto $]\!]J_k^*,\Sigma
_{k+1}^*]\!]$ of the restriction of $\mu$ to the subtree rooted at
$J_k^*$, so that $([\![J_k^*,\Sigma_{k+1}^*]\!],\mu_{k,k+1}^*)$ is a\vspace*{1pt}
rescaled copy of $(\mathcal{R}_1^*,\mu_1^*)$. Since $J_k^*$ is picked
from $\mu_k^*$, we say that $((\mathcal{R}_k^*,\mu_k^*),k\ge1)$
develops by \emph{size-biased branching}, generalising the case $k=1$
that relates to Proposition~\ref{symmbif} via the self-similar time
change (\ref{afreq}).

\subsection{Size-biased bead selection and strongly sampling consistent compositions}\label{sectbsp}

The basic building block for the tree growth process $((\mathcal
{R}_k^*,\mu_k^*),k\ge1)$ is a family of independent copies of
$(\mathcal{R}_1^*,\mu_1^*)$, or equivalently, a family of independent
copies of a fragmenter $M^*$, related by the following
general construction.

%

%
\begin{defn}[(String of beads)]
Given a decreasing pure jump process
$M$ and two positive real parameters $\alpha$ and $m$, we construct a
\emph{string of beads of mass $m$} by placing a random discrete measure
$\mu_{M,\alpha,m}$ on the interval
$(0,Y_{M,\alpha,m}]$ of random length
\[
Y_{M,\alpha,m} = m^\alpha\int_0^\infty
M_s^\alpha \,ds
\]
according to the formula
\[
\mu_{M,\alpha,m} \biggl( m^\alpha\int_0^t
M_s^\alpha \,ds, m^\alpha\int_0^\infty
M_s^\alpha \,ds \biggr) = m M_t. %
\]
If $M$ is a fragmenter with Laplace exponent $\Phi$, we call
$([0,Y_{M,\alpha,m}],\mu_{M,\alpha,m})$ an \emph{$(\alpha,\Phi)$-string of beads of mass $m$}.
\end{defn}

Note that for each $t$ that is a jump time of $M$, the measure $\mu
_{M,\alpha,m}$ puts mass $m(M_{t-} - M_t)$ at the location
$m^\alpha\int_0^t M_s^\alpha \,ds$.
Now, by repeated application of this scheme, we
construct an increasing sequence of $\mathbb{R}$-trees $(R_n,n\ge1)$,
where each $R_n$ is equipped with a random discrete distribution $\mu_n$.

%
\begin{defn}[(Bead splitting process)]\label{defbsp}
Let $\alpha>0$ and $M_n$, $n\ge1$, be a sequence of decreasing pure jump processes
starting from 1:
\begin{itemize}
\item Let $(R_1,\mu_1)$ be the string of beads of mass 1 associated
with $M_1$ and $\alpha$. More specifically, let $R_1=[0,\Sigma
_1]:=[0, Y_{M_1,\alpha,1}]$ be equipped with the usual distance, with
root vertex $0$ and with the random discrete distribution $\mu_1=\mu
_{M_1,\alpha,1}$.
\item Given that $R_n$ has been defined as an $\mathbb{R}$-tree with
root vertex $0$ and $n$ leaves $\Sigma_1,\ldots,\Sigma_n$, and
equipped with a mass measure
$\mu_n$ with total mass $1$, let $R_{n+1}$ be defined as follows. Pick
a junction point $J_n$ from $R_n$ according to
$\mu_n$. Given $\mu_n(\{J_n\}) = m$, distribute the mass $m$
according to a copy
$([\![J_n,\Sigma_{n+1}]\!],\mu_{n,n+1})$ of the string of beads
$([0,Y_{M_{n+1},\alpha,m}],\mu_{M_{n+1},\alpha,m})$, and then attach
this segment to $(R_n,\mu_n-m\delta_{J_n})$ at $J_n$ to form
$(R_{n+1},\mu_{n+1})$.
\end{itemize}
We refer to the projective sequence $((R_n,\mu_n),n\ge1)$ of weighted
$\mathbb{R}$-trees as a \emph{bead splitting process} that
develops by \emph{size-biased branching}.
\end{defn}


We think of $(R_n,\mu_n)$ as $n$ pieces of string $[0,Y_{M_1,\alpha,1}]$, $[0,Y_{M_2,\alpha,\mu(\{J_1\})}],\ldots,\break [0,Y_{M_n,\alpha,\mu
(\{J_{n-1}\})}]$ tied at the junction points $J_1,\ldots,J_{n-1}$,
with beads according to $\mu_n$. The $n$th growth step selects bead
$J_n$ of size $\mu_n(\{J_n\})$ and splits it into smaller beads that
are placed onto a new piece of string tied to $J_n$.

The growth process $((R_n,\mu_n),n\ge1)$ gives rise to an ordered
spinal partition of $\mathbb{N}\setminus\{1\}$ in the terminology of
\cite{HPW}, which we can represent by a point process
%
\begin{equation}
\label{ordsc}\Pi^{\mathrm{ord}}_s=\bigl\{n\in\mathbb{N}\setminus
\{1\}\dvtx  J_{1,n}=g_{0,\Sigma_1}(s)\bigr\},\qquad s\ge0,
\end{equation}
where $J_{1,n}$ is the branch point that has $\Sigma_1$ and $\Sigma
_n$ in two different subtrees, where $[\![0,J_{1,n}]\!]=[\![0,\Sigma_1]\!]\cap
[\![0,\Sigma_n]\!]$, and $g_{0,\Sigma_1}\dvtx [0,d(0,\Sigma
_1)]\rightarrow[\![0,\Sigma_1]\!]$ is the unique isometry with
$g_{0,\Sigma_1}(0)=0$.

%
\begin{prop}\label{regcomp} Given any bead splitting process
$((R_n,\mu_n),n\ge1)$ that develops by size-biased branching, the
ordered spinal
partition $\Pi^{\mathrm{ord}}$ defined in (\ref{ordsc}) is exchangeable. In
particular, if we
choose $M_1$ to be a fragmenter, then the spinal partition gives rise
to a (strongly sampling consistent) regenerative
composition structure $(\mathcal{C}_n,n\ge1)$, which records for each
$n\ge1$ the vector $\mathcal{C}_n$ of nonzero block sizes
$\#(\Pi^{\mathrm{ord}}_s\cap[n+1])$, $s\ge0$, of $\Pi^{\mathrm{ord}}\cap
[n+1]$, arranged in the spinal order of blocks given by the
order of positions $s\ge0$ on the spine.
\end{prop}

\begin{pf} The first statement holds since $\mu_1$ is the projection
of $\mu_n$ to $R_1$ for all $n\ge1$, so the picks of
$J_n$ projected to $R_1$ are exchangeable picks from $\mu_1$ by the
use of size-biased branching. The second statement now follows directly
from Gnedin and Pitman
\cite{gp}, Theorem 5.2; cf. also \cite{PW09}, Section~2.1, for the
terminology of (weak and) strong sampling consistency.
\end{pf}

\subsection{Convergence of bead splitting processes to self-similar CRTs}

The next theorem establishes CRT convergence of bead splitting
processes\break $((R_n,\mu_n),n\ge1)$ in the sense of Definition~\ref
{defbsp}, not just for the case of symmetric splitting rules $f^*$ that
relate directly to the growth process $((\mathcal{R}_n^*,\mu
_n^*),n\ge1)$ obtained by sampling from the measure $\mu$ of a CRT
$(\mathcal{T},\mu)$, but also for fragmenters $(M_n,n\ge1)$ with
nonsymmetric splitting rules $f$, with convergence to a CRT associated
with the symmetrised splitting rule $f^*$ associated with $f$. Again,
the result holds without assuming the existence of densities $f$ and $f^*$:

%
\begin{teo}\label{teogrowth}
For independent fragmenters $M_n$ with Laplace exponent $\Phi(\rho
)=\int_{(0,1)}(1-u^\rho)\Lambda(du)$ and $\alpha>0$, the sequence
of weighted random $\mathbb{R}$-trees $(R_n,\mu_n)$ converges
almost surely in the Gromov--Hausdorff--Prohorov metric to a limit tree
$(\mathcal{T},\mu)$, which is a copy of the $\alpha$-self-similar
tree that is canonically
associated with a binary fragmentation process with symmetric
dislocation measure $\nu=\Lambda+\overline{\Lambda}$. In addition,
we also have
$(R_n,\nu_n)\rightarrow(T,\mu)$ almost surely in the
Gromov--Hausdorff--Prohorov metric, where $\nu_n$ is the uniform
measure on the
$n$ leaves of $R_n$.
\end{teo}

If\vspace*{1pt} $f(u) =\frac{1}{\sqrt{2\pi}} u^{-3/2}(1-u)^{-3/2}$, it follows
easily from the work of Haas and Miermont \cite{HM04} that the
sequence $(R_n,n\ge1)$
has the same distribution as the increasing sequence of trees provided
by Aldous's \cite{Ald-crt1} line-breaking construction of the Brownian
CRT. Therefore, Theorem~\ref{teogrowth} can be seen as a
generalisation of Aldous's line-breaking construction. We discuss this
example of a bead splitting process $((R_n,\mu_n),n\ge1)$ in
Section~\ref{bcrt}.

To prove this theorem, we will embed $(R_n,\mu_n)$ in a CRT $(T,\mu
)$, as has (essentially) been done for $(R_1,\mu_1)$ in Theorem~\ref
{embed}. A key tool will be the following spinal decomposition result.

%
\begin{lm}[(Spinal decomposition)]\label{spindecA} Let $A=(A_t,t\ge0)$
be a Markovian path in a homogeneous fragmentation process
$\Pi=(\Pi(t),t\ge0)$. For each
$n\ge1$, denote by $\Pi^{\{n\}}(t)$ the block of $\Pi(t)$ containing
$n$, $t\ge0$, consider
$\sigma_n=\inf\{t\ge0\dvtx  n\notin A_t\}$ and the associated spinal
partition $\Pi^A(0)=\{\Pi^{\{n\}}(\sigma_n),n\ge1\}$. Then
conditionally given $\Pi^A(0)$ and $(\sigma_n,n\ge1)$, the
process
\[
\Pi^A(t)=\bigl\{\Pi^{\{n\}}(\sigma_n+t),n\ge1
\bigr\}=\bigcup_{i\ge1}\Pi(\sigma_{\min\Pi^A_i(0)}+t)
\cap\Pi_i^A(0),\qquad t\ge0
\]
is a fragmentation process starting from $\Pi^A(0)$, with the same
transition kernel as~$\Pi$.
\end{lm}

This lemma says that the process $\Pi$ can be decomposed along the
path $A$ into the partition $\Pi^A(0)$ of blocks that separate from
$A$ at any time $t\ge0$. The blocks $\Pi_i^A(0)$, $i\ge1$ then
evolve independently and according the transition kernel of $\Pi$.

\begin{pf*}{Proof of Lemma~\ref{spindecA}}
We extend the proof of \cite{HPW}, Proposition~4, to the
higher generality here of decomposing along a Markovian path.
The family of times $(\sigma_n,n\ge1)$ is a stopping line for the
filtration $\mathcal{F}=(\mathcal{F}_t,t\ge0)$, with respect to
which $A$ is a
Markovian path in $\Pi$. Here, we use the terminology of Bertoin \cite
{ber-book}, Definition 3.4, and seek to obtain from
\cite{ber-book}, Lemma 3.14, that the \emph{extended branching property} holds, which yields precisely the result we need. Since
Bertoin uses natural filtrations, and to demonstrate where the
Markovian assumption on the path enters the argument, let us
briefly retrace Bertoin's steps and sketch relevant parts of the proof
of the extended branching property. Without loss of
generality, $\mathcal{F}$ is the filtration generated by $(\Pi,A)$.
Also denote by $\mathcal{F}^{\{n\}}$ the filtrations generated
by $(\Pi^{\{n\}},A_{ \cdot\wedge\sigma_n})$, for each $n\ge1$.
We\vspace*{1pt} consider approximations $\sigma_n^{(h)} =\inf\{t\in h\mathbb
{N}\dvtx  n\notin A_t\}$, $\sigma_n^{(h,k)} =\min(kh,\sigma
_n^{(h)})$, and
$\overline{\sigma}_n^{(h,k)} =\inf\{t\in\{h,2h,\ldots,kh\}\dvtx
n\notin A_t\}\}$ with $\inf\varnothing:=\infty$, of $\sigma_n$.
The branching property at the stopping line $(\sigma_n^{(h,1)},n\ge
1)$ is just the branching property at $t=h$. At $h$, or
by induction hypothesis at $(\sigma_n^{(h,k)},n\ge1)$, the assumption
on the path to be Markovian ensures that
$(\Pi(kh+t)\cap A_{kh},A_{kh+t},t\ge0)$ is conditionally independent of
$(\{\Pi^{\{n\}}(\sigma^{(h,k)}_n+t),n\notin A_{kh}\},t\ge0)$ given
\[
\mathcal{F}_{(\sigma^{(k,h)}_n,n\ge1)},\mbox{ defined as the
sigma-algebra generated by }
\mathcal{F}_{\sigma_n^{(h,k)}}^{\{n\}}, n\ge1.
\]
To
$(\Pi_{kh+t}\cap A_{kh},A_{kh+t},t\ge0)$, we can apply the branching
property at $t=h$ and trivially at $t=\infty$ to complete
the induction step from $k$ to $k+1$. This establishes the extended
branching property at $(\sigma_n^{(h,k)},n\ge1)$ and
$(\overline{\sigma}_n^{(h,k)},n\ge1)$ for all $h>0$ and $k\ge1$. We
omit the remainder of the proof, which
uses the standard approximation $\overline{\sigma
}_n^{2^{-k},2^{2k}}\downarrow\sigma_n$ as $k\rightarrow\infty$.
\end{pf*}

The next lemma and its proof demonstrate that we can iterate the
embedding of a Markovian path in a homogeneous fragmentation process to
embed a bead splitting process in an associated self-similar CRT to
which the bead splitting process converges almost surely.

%
\begin{lm}\label{markpathconv} Let $A=(A_t,t\ge0)$ be a Markovian
path in a binary fragmentation process $\Pi=(\Pi(t),t\ge0)$ and
$M_t=|A_t|$, $t\ge0$ its residual mass process. If the $M_n$, $n\ge1$ are
independent copies of $M$, then for each $\alpha> 0$ the sequence
of weighted random $\mathbb{R}$-trees $(R_n,\mu_n)$ converges almost
surely in the Gromov--Hausdorff--Prohorov metric to a limit tree
$(\mathcal{T},\mu)$, which is a copy of the $\alpha$-self-similar
CRT that
is canonically associated with $\Pi$. In addition, we also have
$(R_n,\nu_n)\rightarrow(\mathcal{T},\mu)$ almost surely in the
Gromov--Hausdorff--Prohorov metric, where $\nu_n$ is the uniform
measure on
the $n$ leaves of $R_n$.
\end{lm}

\begin{pf} We can consider the spinal partition $\Pi^A(0)$ of $\Pi$
and use Lemma \ref{spindecA} to construct in a measurable way
(see, e.g., \cite{HPW}, Corollary~3) a string of beads $(R_1,\mu
_1)=([0,Y_{|A|,\alpha,1}],\mu_{|A|,\alpha,1})$ with a
collection $(s_i,T^{(i)},\mu^{(i)})$ of spinal subtrees constructed from
$(\Pi(\sigma_{\min\Pi^A_i(0)}+t)\cap\Pi^A_i(0),t\ge0)$, where\vspace*{-3pt}
$s_i=\int_0^{\sigma_{\min\Pi^A_i(0)}}|A_r|^\alpha \,dr$, $i\ge1$,
such that the tree
$(\mathcal{T},\mu)$ obtained by grafting $(T^{(i)},\mu^{(i)})$ to
$(R_1,0)$ at $s_i$ for all $i\ge1$, is a self-similar CRT associated
with $\Pi$. This construction gives rise to a family of regular
conditional distributions of
$(R_1,\mu_1;(s_i,T^{(i)},\mu^{(i)}),i\ge1)$ given $(\mathcal{T},\mu
)$, and we can
use these via the Ionescu--Tulcea theorem to obtain a probability space
that allows the following construction.

Suppose we have constructed $(R_n,\mu_n;(x_i,T^{(i)},\mu^{(i)}),i\in
I_n)$ with $R_n\subset\mathcal{T}$, $\mu_n$ the projection of $\mu$
onto $R_n$ and, conditionally given $(R_n,\mu_n)$, a collection
$((T^{(i)},\mu^{(i)}),i\in I_n)$ of independent copies of
$(\mathcal{T},\mu)$ subjected to $\alpha$-self-similar scaling by
$\mu(\{x_i\})$, which when grafted at $x_i\in R_n$ for all $i\in I_n$
give $(\mathcal{T},\mu)$. Now pick a junction point $J_n=x_{i_n}$
from $R_n$ according to $\mu_n$. Given that $\mu_n(\{J_n\})=m$, remove
$J_n$ from $\mu_n$ and remove $i_n$ from $I_n$. Use the regular
conditional distribution given the rescaled chosen subtree
$(T^{(i_n)},\mu^{(i_n)})$ to obtain a string of beads with grafted
spinal subtrees distributed as
$(R_1,\mu_1;(s_i,T^{(i)},\mu^{(i)}),i\ge1)$, without
modifying the chosen rescaled subtree. After $\alpha$-self-similar
scaling by $m$, graft the string of beads at $J_n$, add the
new spinal subtrees to the collection to form $(R_{n+1},\mu
_{n+1};\break (x_i,T^{(i)},\mu^{(i)}),i\in I_{n+1})$.
Then $R_{n+1}\subset\mathcal{T}$, $\mu_{n+1}$ is the projection of
$\mu$ onto $R_{n+1}$ and, conditionally given $(R_{n+1},\mu_{n+1})$,
the collection $((T^{(i)},\mu^{(i)}),i\in I_{n+1})$ consists of scaled
independent copies of $(\mathcal{T},\mu)$ that turn $R_{n+1}$ into
$(\mathcal{T},\mu)$ when grafted at $x_i$, $i\in I_{n+1}$.

By induction, this gives a sequence $((R_n,\mu_n),n\ge1)$ embedded in
$(\mathcal{T},\mu)$, which develops by size-biased branching and
is based on independent copies of the string of beads associated with
$M=|A|$. While constructed within a CRT, this sequence has
the same (joint) distribution as the sequence described in the
statement of the lemma. It therefore suffices to prove almost
sure convergence for this embedded sequence.

First, consider the measures $\mu_n$, $n\ge1$. Denote by $|\mu
_n|^\downarrow\in\mathcal{S}^\downarrow$ the decreasing
rearrangement of
bead sizes $\mu_n(\{x\})$, $x\in R_n$. Since $A$ is embedded in
$(\mathcal{T},\mu)$, the measure $\mu_1$ cannot have an atom of size 1;
in particular there is $\lambda<1$ such that $\mathbb{P}(|\mu
_1|^\downarrow_1<\lambda)>0$. Now let $\varepsilon=1/K>0$. By selecting
the $K$ largest beads in turn, we see that $\mathbb{P}(|\mu
_{n+K}|^\downarrow_1<\lambda s_1\vee\varepsilon\giv|\mu
_n|^\downarrow=\mathbf{s})>0$
for all $\mathbf{s}\in\mathcal{S}^\downarrow$. For $m$ with
$\lambda
^m<\varepsilon$ this implies $p=\mathbb{P}(|\mu_m|^\downarrow
_1<\varepsilon)>0$,
but then $|\mu_n|^\downarrow_1$ will be less than $\varepsilon$
after a time that is bounded above by $m$ times a geometric
random variable with parameter $p$. In particular,
%
\begin{equation}
\label{atom}\mathbb{P}(\mu_n\mbox{ has an atom of size greater than
$\varepsilon$ for all $n\ge1$})=0.
\end{equation}
Now denote by $R_\infty$ the completion of the increasing union
$\bigcup_{n\ge1}R_n$ in $\mathcal{T}$, and assume that
$\mathbb{P}(R_\infty\neq\mathcal{T})>0$. For $x\in\mathcal
{T}\setminus R_\infty$, we find $y\in R_\infty$ such that
$]\!]y,x]\!]\cap R_\infty=\varnothing$,
but then $\mu(\mathcal{T}_y)>0$, since $\mu$ assigns positive weight
to all fringe subtrees. Since $\mu_n$ is the projection of $\mu$
onto $R_n\subset R_\infty$, this contradicts (\ref{atom}). Hence
$\mathbb{P}(R_\infty=\mathcal{T})=1$. Similarly, assuming that $R_n$
does not
converge to $\mathcal{T}$ for the Hausdorff distance on $\mathcal
{T}$, we can use compactness to find $x\in\mathcal{T}$ with
$d(x,R_n)>\varepsilon$ for
all $n\ge1$, so $x\notin R_\infty$ is a contradiction. Also,
Hausdorff convergence of $R_n\subset\mathcal{T}$ to $\mathcal{T}$ with
projected measures implies Gromov--Hausdorff--Prohorov convergence
$d_{\mathrm{GHP}}(R_n,\mathcal{T})\rightarrow0$ almost surely as
$n\rightarrow\infty$; see, for \mbox{example}, \cite{PW09}, Lemma 17.

Finally, the measure $\nu_n$ on the $n$ leaves on $R_n$ is more and
more closely coupled with the measure $\nu_n^*$ on $\mathcal{T}$
associated with a sample $\Sigma_1^*,\ldots,\Sigma_n^*$ from $\mu$,
which is well known to
converge weakly almost surely to $\mu$. We can take as $\Sigma_1^*$
an independent pick from $\mu$
and include $\Sigma_{n+1}^*$ in the construction of $(R_{n+1},\mu
_{n+1})$. Specifically, we can
obtain the pick from $\mu_n$ for the junction point $J_n$ as the
junction point of the subtree containing $\Sigma_{n+1}^*$. Since
$R_n\rightarrow\mathcal{T}$ almost surely in the Hausdorff sense,
there is $n_0$ such that all subtrees of $\mathcal{T}\setminus
R_{n_0}$ have
height less than $\varepsilon/2$, but then the distance between
$\Sigma_{n+1}^*$ and the $(n+1)$st leaf of $R_{n+1}$, which are
in the same subtree by construction, is at most $\varepsilon$ for all
$n\ge n_0$, which entails the result by standard arguments.
\end{pf}

\begin{pf*}{Proof of Theorem~\ref{teogrowth}}
Let $M$ be a fragmenter with Laplace exponent $\Phi(\rho)=\int
_{(0,1)}(1-u^\rho)\Lambda(du)$. By the
existence part of Theorem~\ref{embed}, there is a Markovian embedding
of $M$ into $\Pi$ for a binary homogeneous fragmentation
process with symmetric dislocation measure $\nu=\Lambda+\overline
{\Lambda}$. Hence, Lemma~\ref{markpathconv} applies and gives
$(R_n,\mu_n)\rightarrow(\mathcal{T},\mu)$ and $(R_n,\nu
_n)\rightarrow(\mathcal{T},\mu)$ almost surely in the
Gromov--Hausdorff--Prohorov sense,
for an $\alpha$-self-similar CRT $(\mathcal{T},\mu)$ with symmetric
dislocation measure $\nu$.
\end{pf*}

This establishes Theorem~\ref{teogrowth}, and in particular shows that
the bead-splitting processes of Definition~\ref{defbsp} answer points
4--5 from the \hyperref[sec1]{Introduction}, and that point~6 then also holds. We
can now use Theorem~\ref{teogrowth} and Lemma~\ref{markpathconv} to
complete the proof of Theorem~\ref{embed}, and hence establish point
2, completing the programme of 6 points set out in the \hyperref[sec1]{Introduction}:

\begin{pf*}{Proof of the uniqueness part of Theorem~\ref{embed}}
Let\vspace*{1pt} $\Phi(\rho)=\break\int_{(0,1)}(1-u^\rho)\Lambda(du)$ be the Laplace
exponent of a fragmenter. Consider the bead
splitting process $((R_n,\mu_n),n\ge1)$ based on a sequence $M_n$,
$n\ge1$, of such fragmenters. In Theorem~\ref{teogrowth} we showed that
$(R_n,\mu_n)\rightarrow(\mathcal{T},\mu)$ almost surely for a CRT
$(\mathcal{T},\mu)$ with symmetrised dislocation measure $\nu
=\Lambda+\overline{\Lambda}$.

Now assume that a fragmenter with Laplace exponent $\Phi$ has a
Markovian embedding $A$ into an exchangeable binary homogeneous
fragmentation process $\Pi$ with any symmetric dislocation measure
$\widetilde{\nu}$. By Lemma~\ref{markpathconv},
$(R_n,\mu_n)\rightarrow(\widetilde{\mathcal{T}},\widetilde{\mu})$
almost surely for a CRT with symmetrised dislocation measure
$\widetilde{\nu}$. By uniqueness of limits, $(\widetilde{\mathcal
{T}},\widetilde{\mu})=(\mathcal{T},\mu)$. Since the distributions
of CRTs for
different dislocation measures are different, we find that $\nu
=\widetilde{\nu}$.
\end{pf*}

In the remainder of the paper, we point out some further connections to
related work.

%
\begin{rem}\label{nexhom}
With the usual names $\Sigma_1,\ldots,\Sigma_n$ of leaves of $R_n$, $n\ge1$, any bead splitting process
embedded in a CRT $(\mathcal{T},\mu)$ gives rise to a, typically
nonexchangeable, $\mathcal{P}$-valued process
\[
\Pi_\alpha(t)=\bigl\{\bigl\{j\ge1\dvtx \Sigma_j\in
\mathcal{T}^t_i\bigr\},i\ge1\bigr\}\cup\bigl\{\{j\},j
\ge1\dvtx \Sigma_j\notin\mathcal{T}^t\bigr\},\qquad t
\ge0, %
\]
of the same form as the exchangeable special case $\Pi_\alpha^*$ in
(\ref{embedpart}). Furthermore, if $(R_n,\mu_n)\rightarrow(\mathcal
{T},\mu)$ as in Lemma~\ref{markpathconv}, we find equality of the
decreasing rearrangements of asymptotic frequencies $|\Pi_\alpha
(t)|^\downarrow=|(\mu(\mathcal{T}^t_i),i\ge1)|^\downarrow=|\Pi
_\alpha^*(t)|^\downarrow$ for all $t\ge0$ a.s. So it is natural to
perform the inverse of the self-similar time change (\ref{afreq}) to
construct a, typically nonexchangeable, homogeneous process $\Pi=\Pi
_0$ from the consistently time-changed evolution of its blocks
containing $n$, $n\ge1$.
\end{rem}

Our proof of the uniqueness part of Theorem~\ref{embed} used the
size-biased bead splitting process and the compactness of self-similar
CRTs to show that the Markovian path $A$ gives rise to an embedding
that exhausts a CRT. The embedding for the existence part of Theorem
\ref{embed} was not carried out in a CRT, but directly in an
exchangeable homogeneous fragmentation process. Indeed, it should be
possible to also prove the uniqueness in the framework of homogeneous
fragmentation processes. We can rephrase our bead splitting argument
for the uniqueness part here to directly construct a nonexchangeable
process $\Pi$ based on $A$ by embedding into an exchangeable
homogeneous fragmentation process $\Pi^*$, as indicated below.
However, this is harder to formulate, and we
lose natural compactness, so we do not attempt an alternative proof,
but let us give the direct construction of $\Pi$.

Let $A$ be a Markovian path in $\Pi^*$. Define branch times
$J_{1,n}=\inf\{t\ge0\dvtx  n\notin A_t\}$, $n\ge2$, between $1$ and
$n$. Given $J_{i,n}$, $n\ge i+1$, have
been constructed for all $i\in[k]:=\{1,\ldots,k\}$, consider the time
$H_{k+1}=\max\{J_{i,k+1},i\in[k]\}$ when $k+1$ separates from the
last $i\in[k]$. Relabel the restriction of $(\Pi^*(H_{k+1}+t),t\ge
0)$ to the block $B_{k+1}$ of $\Pi^*(H_{k+1})$ that contains $k+1$ by
the increasing bijection $B_{k+1}\rightarrow\mathbb{N}$. Run a copy
of $A$ inside this process,\vspace*{1pt} relabel back $\mathbb{N}\rightarrow
B_{k+1}$ to find a Markovian path $B^{(k+1)}$ that we specify to
coincide with the canonical path $A^{k+1}$ of $\Pi^*$ up to $H_{k+1}$
and to
continue in $(\Pi^*(H_{k+1}+t)\cap B_{k+1},t\ge0)$, as constructed.
Define $J_{k+1,n}=\inf\{t\ge0\dvtx  n\notin B^{(k+1)}_t\}$, $n \ge
k+2$. Finally set $J_{k,k}=\infty$, $J_{k,n}=J_{n,k}$ for $n<k$ and
define the embedded $\mathcal{P}$-valued process
%
\begin{equation}
\label{nonexpi}\Pi_t=\bigl\{\{n\ge1\dvtx  J_{n,i}>t\},i\ge1
\bigr\},\qquad t\ge0.
\end{equation}

%
\begin{cor} Let $A$ be a Markovian path in an exchangeable homogeneous
fragmentation process $\Pi^*$, and let $\Pi$ be as in (\ref{nonexpi}).
Then $|\Pi|^\downarrow=|\Pi^*|^\downarrow$. Moreover, if $|A|$ is a
fragmenter, then $\Pi$ is a homogeneous fragmentation
process with binary nonexchangeable $\kappa$-measure \cite
{PRW} of the form
\begin{eqnarray}
\kappa\bigl(\bigl\{\Gamma\in\mathcal{P}\dvtx \Gamma\cap[n]=(\pi_1,
\pi_2)\bigr\}\bigr)=\int_{(0,1)}u^{\#\pi_1-1}(1-u)^{\#\pi_2}
\Lambda(du),\nonumber
\\
\eqntext{\{ \pi_1,\pi_2\}\in\mathcal{P}_n
\setminus\bigl\{\bigl\{[n]\bigr\}\bigr\},}
\end{eqnarray}
where $\mathcal{P}_n$ is the set of partitions of $[n]:=\{1,\ldots,n\}
$, with $\kappa(\Gamma\in\mathcal{P}\setminus\{\{\mathbb{N}\}\}
\dvtx \Gamma_1\cup\Gamma_2\neq\mathbb{N})=0$.
\end{cor}

\begin{pf} We leave the equivalence of the two constructions of $\Pi$
to the reader and just point out that the CRT construction
of Remark~\ref{nexhom} yields $|\Pi|=|\Pi^*|$. For the second claim,
we note that the fragmenter has Laplace exponent $\Phi(\rho)=\int
_{(0,1)}(1-u^\rho)\Lambda(du)$, so, by
standard thinning properties of the Poisson point process of jumps of
the fragmenter $M$ and size-biased branching, we identify
the dislocation measure.
\end{pf}

%
\begin{rem}
It may be observed from the form of the bead splitting
process in the case of an independent and identically distributed
sequence $(M_n,n\ge1)$ that the size-biased bead selection rule is not
crucial for convergence to a CRT since it mainly affects the (random)
time $n$ at which a particular bead is split. In the proof of Lemma
\ref{markpathconv}, the main use of the size-biased selection rule was
to establish (\ref{atom}). Indeed, as long as we split every bead
eventually, we are quite free to choose the order in which we split the
beads and may even contemplate rules like splitting all beads of $\mu
_n$ at once at every stage of the bead splitting process.

The reader may also want to compare our bead splitting processes with
Abraham's \cite{Abr-92} construction of a version of the Brownian CRT.
Let us rephrase Abraham's construction in our present framework. The
construction is based on the distribution of the total height of the
CRT and the decomposition of the CRT along the path from the root to
the \emph{highest} leaf. If we project the mass measure of the CRT
onto the spine, we obtain a string of beads $(\overline{R}_1,\overline
{\mu}_1)$. Abraham takes this string of beads and recursively replaces
all beads of $\overline{\mu}_n$ by a rescaled copy of $(\overline
{R}_1,\overline{\mu}_1)$ conditioned not to exceed the height of the
branch of the bead. The path to the highest leaf does not correspond to
a Markovian path, and the spinal subtrees are not rescaled copies of
the CRT, but copies constrained in height, this falls outside the
setting of Lemma~\ref{markpathconv}.

Intuitively, the \emph{Markovian} path $A$, whose leaf height
$L_{0\Sigma}=\int_0^\infty|A_t|^\alpha \,dt$ in an associated $\alpha
$-self-similar CRT is most likely to be highest is the one based on the
switching probabilities of Example~\ref{ex7}, always choosing the
bigger fragment. The homogeneous Poissonian structure for relative
masses $(F_t,1-F_t)$ on the spine to the embedded leaf easily entails
that this does not lead to the highest leaf a.s.
\end{rem}

\subsection{The Brownian CRT}\label{bcrt}
Let $(\mathcal{T},\mu)$ denote the Brownian Continuum Random Tree
$\mathcal{T}$ equipped with its mass measure $\mu$,
which Aldous \cite{Ald-crt3} constructed both as the tree embedded in
(twice the standard)
Brownian excursion, with $\mu$ corresponding to Lebesgue measure on $[0,1]$,
and as a limit as $n \rightarrow\infty$ of an increasing sequence of binary
subtrees with edge-lengths $R_n$
with $n$ leaves labelled by $[n]$, in which case $\mu$ may be interpreted
as the almost sure weak limit as $n\rightarrow\infty$ of the uniform
probability distribution
$\nu_n$ on the $n$ leaves of $R_n$.
The tree $R_n$ may be constructed as the subtree of $\mathcal{T}$ spanned
by $n$ leaves of $\mathcal{T}$, which given $\mathcal{T}$ are picked
independently
according to the mass measure $\mu$. We recover this second
construction in Theorem~\ref{teogrowth} for $f(u)=\frac{1}{\sqrt
{2\pi}}u^{-3/2}(1-u)^{-3/2}$, enriched by the string of beads
structure given to the branches of $R_n$ by measures $\mu_n$.
According to a basic result of Aldous \cite{Ald-crt3}, the increasing
sequence of lengths $(\lambda(R_1),\lambda(R_2), \ldots)$ of these
subtrees can be constructed as $\lambda(R_n) = \sqrt{2 \Gamma_n }$
where $\Gamma_n = \varepsilon_1 + \cdots+\varepsilon_n $ for a
sequence of independent
standard exponential variables $\varepsilon_i$.
For $n=1$ there is the identity in distribution
%
\begin{equation}
\label{rayleigh} \lambda(R_1) \stackrel{d} {=}L^0_1
\bigl(B^{\mathrm{br}}\bigr),
\end{equation}
where $B^{\mathrm{br}}$ is a standard Brownian bridge, starting at $0$
at time
$0$ and ending at~$0$ at time $1$, and
$(L_t^x(B^{\mathrm{br}}), 0 \le t \le1, x \in\mathbb{R})$ is the jointly
continuous process
of local times of $B^{\mathrm{br}}$, normalised so that
$L_t^x(B^{\mathrm{br}})\,dx$ is the
occupation measure of $(B^{\mathrm{br}}(s),0\le s\le t)$. The common
distribution
of both sides
in (\ref{rayleigh}) has the Rayleigh density $x\exp(-x^2/2)$ at $x
\ge0$.
The work of Aldous, Miermont, and Pitman \cite{AMP-04a} yields
a deeper connection between the Brownian CRT in a Brownian excursion on
the one hand and
Brownian bridge on the other. This establishes a spinal decomposition
result for the Brownian CRT via path transformations rather than via
Bertoin's extended branching property as in \cite{HPW} or Lemma~\ref
{spindecA} here. See also Bertoin and Pitman \cite{BeP-94},
Theorem~3.2, for an expression in terms of paths rather than trees.
To explore this spinal decomposition, let $(\tau_\ell, 0 \le\ell<
L^0_1(B^{\mathrm{br}}))$ be the inverse local time process
\[
\tau_\ell:= \inf\bigl\{t\ge0\dvtx  L_t^0
\bigl(B^{\mathrm{br}}\bigr) > \ell\bigr\} %
\]
so that the collection of excursion intervals of $B^{\mathrm{br}}$ is
exhausted by
\[
\bigl\{ (\tau_{\ell- }, \tau_\ell)\dvtx  \ell>0,
\tau_{\ell- } < \tau_\ell\bigr\}, %
\]
and let $P$ be the random discrete distribution obtained by ranking these
intervals by length.
On the other hand, in the Brownian CRT $(\mathcal{T},\mu)$,
for $0 \le\ell< \lambda(R_1)$ let
$\mu_1([0,\ell])$ denote the mass of all points $x \in\mathcal{T}$
such that the
path from root to~$x$ in $\mathcal{T}$ branches from the path from
root to leaf
$\Sigma_1$ of~$\mathcal{T}$ at a junction point on~$R_1$ whose
distance from the root
of $\mathcal{T}$ is at most $\ell$.
Then according to the spinal decomposition of the Brownian CRT
implied by \cite{AMP-04a}, Lemma 9 and equation (12), the equality in
distribution (\ref{rayleigh})
extends to the equality in distribution of processes
%
\begin{equation}
\label{locid} \bigl(\mu_1\bigl([0,\ell]\bigr), 0 \le\ell<
\lambda(R_1) \bigr) \stackrel{d} {=} \bigl(\tau_\ell, 0 \le
\ell< L^0_1\bigl(B^{\mathrm{br}}\bigr)\bigr).
\end{equation}
Moreover, conditionally given the process on the left-hand side
of (\ref{locid}), the Brownian CRT $(\mathcal{T}, \mu)$ decomposes into
a countable collection of subtrees
\[
\bigl\{ \mathcal{T}_\ell, 0 < \ell< \lambda(R_1),
\mu_1\bigl(\{\ell\}\bigr)>0\bigr\},
\]
where $\mathcal{T}_\ell$ is a Brownian CRT equipped with a mass
measure $\mu_\ell$
with total mass $\mu_\ell(\mathcal{T}_\ell) = \mu_1(\{\ell\})>0$.
This decomposition corresponds on the right-hand side of~(\ref{locid}) to
(trees in) excursions of $|B^{\mathrm{br}}|$, the absolute value of
$B^{\mathrm{br}}$, excursions
of lengths
$\{ \tau_{\ell} - \tau_{\ell-}, 0 < \ell< L^1_0(B^{\mathrm{br}}),
\tau_{\ell
- } < \tau_\ell\}$.
Indeed, the entire Brownian CRT $\mathcal{T}$ can be constructed from
$B^{\mathrm{br}}$
so that the equality in distribution (\ref{locid}) holds
almost surely, and for each $\ell$ with
$\mu_1(\{\ell\}) = \tau_\ell- \tau_{\ell-} > 0$,
the subtree $\mathcal{T}_\ell$ of $\mathcal{T}$ attached to the spine
$R_1$ of $\mathcal{T}$ at distance $\ell$ from the root is
constructed from the
excursion of $|B^{\mathrm{br}}|$ on $(\tau_{\ell-},\tau_\ell)$ in the
same way that the entire tree $\mathcal{T}$ is constructed from a standard
Brownian excursion.
In particular, basic properties of Brownian excursions then imply the
\emph{spinal decomposition of $\mathcal{T}$}, that conditionally given the
subtree masses $(\mu_1(\{\ell\}),0 \le\ell< \lambda(R_1) )$, the
subtrees $(\mathcal{T}_\ell,\mu_\ell)$ associated with $\ell$ such that
$\mu_1(\{\ell\})>0$ form
a collection of independent random trees distributed like
$(\sqrt{\mu_1(\{\ell\})} \mathcal{T}, \mu_1(\{\ell\})\mu)$,
meaning\vspace*{1pt} that all edge-lengths in $\mathcal{T}$ are scaled by a factor of
$\sqrt{\mu_1(\{\ell\})}$, while all masses are scaled by
a factor of $\mu_1(\{\ell\})$.

The distribution of ranked masses of atoms of $\mu_1$ is the
Poisson--Dirichlet distribution
$\mathrm{PD}(\frac{1}{2},\frac{1}{2})$, which is the distribution of
ranked lengths of
excursions of Brownian bridge, and masses
corresponding to these lengths are distributed along the spine of length
$\lambda(R_1)$ in an exchangeable random order. Moreover, the length
$\lambda(R_1)$ is itself a measurable functional of the $\mathrm
{PD}(\frac{1}{2},\frac{1}{2}
)$ random\vspace*{1pt}
discrete distribution of masses along the spine, as discussed
in \cite{PPY,csp}.

For $n\ge1$, the bead splitting process $((R_n,\mu_n),n\ge1)$
described in terms of decreasing mass processes in
Definition~\ref{defbsp} can be described in terms of Brownian bridges,
as follows:

\begin{itemize}
\item Start from a segment
%
\begin{equation}
\label{bsb}(R_1,\mu_1)=\bigl(\bigl[0,L_1^0
\bigl(B^{\mathrm
{br}}\bigr)\bigr],d\tau\bigr),
\end{equation}
where $d\tau$ denotes the Stieltjes measure with cumulative
distribution function given in (\ref{locid}), associated with the inverse
local time $\tau$ of a standard Brownian bridge $B^{\mathrm{br}}$ of
length $1$.

\item Given $(R_n,\mu_n)$, pick a junction point $J_n$ from $R_n$
according to $\mu_n$. Given that $\mu_n(\{J_n\})=m$,
remove the mass $m$ from point $J_n$ and attach as segment
$(]\!]J_n,\Sigma_{n+1}]\!],\mu_{n,n+1})$ a copy of
(\ref{bsb}) derived from a Brownian bridge of length~$m$.
\end{itemize}
Specifically, a Brownian bridge of length $m$ may be constructed from
the standard Brownian bridge $B^{\mathrm{br}}$ as $m^{1/2} B^{\mathrm
{br}}(t/m)$, $0 \le t
\le m$.
That is to say, $R_{n+1}\setminus R_n=]\!]J_n,\Sigma_{n+1}]\!]$ is such that
\[
\lambda(R_{n+1}) - \lambda(R_n) = \bigl(\mu_n
\bigl(\{J_n\}\bigr)\bigr)^{1/2} L_1^0
\bigl( B^{\mathrm{br}} _{(n)} \bigr) %
\]
for some independent and identically distributed sequence of standard
Brownian bridges
\[
B^{\mathrm{br}}_{(n)} = \bigl(B^{\mathrm{br}}_{(n)}(t), 0
\le t \le1 \bigr), %
\]
and given $\mu_n(\{J_n\}) = m$, the mass $m$ should be reallocated with
a portion $(\tau_\ell- \tau_{\ell- })m$ placed at distance
$\ell m^{1/2}$ from $J_n$ along\vspace*{1pt} the branch of length
$\lambda(R_{n+1}) - \lambda(R_n)$ from $J_n$ to $\Sigma_{n+1}$, for each
$\ell\in(0, L_1^0 ( B^{\mathrm{br}}_{(n)} ))$ with $\tau_{\ell-} <
\tau_{\ell}$.\vspace*{1pt}

The above prescription specifies
the projective sequence of
weighted ${\mathbb{R}}$-trees
$( (R_n, \mu_n), n \ge1)$ as in Theorem~\ref{teogrowth}, which is
associated with
the Brownian CRT as in \cite{HM04} and as indicated at the beginning
of Section~\ref{genline}.
We wish to point out some
special properties of this Brownian tree growth sequence:

%
\begin{prop}
Let $((R_n,\mu_n),n\ge1)$ be the bead splitting process derived from
a sequence of Brownian bridges, as above.
Then we have the following description of the law of $(R_n,\mu_n)$:

\begin{longlist}[(iii)]
\item[(i)] The\vspace*{1pt} sequence $P_n = (P_{n,1}, P_{n,2}, \ldots)$ of
sizes of ranked atoms of $\mu_n$ has
$\mathrm{PD}(\frac{1}{2}, (2n-1)/2)$ distribution.\vspace*{1pt}

\item[(ii)] The total length $\lambda(R_n)$ can be represented as both
$\lambda(R_n) = \sqrt{2 \Gamma_n }$ where $\Gamma_n = \varepsilon
_1 + \cdots+\varepsilon_n $
for a sequence of independent standard exponential variables~$\varepsilon_i$,
and as $\lambda(R_n) = S_{1/2} (P_n)/\sqrt{2}$, where $S_{1/2}(P_n)$ is
the $\frac{1}{2}$-diversity of $P_n$, which may be recovered from
$P_n$ as
\[
S_{1/2} (P_n) =\sqrt{\pi}\lim_{k \rightarrow\infty} k
P_{n,k}^{1/2} \mbox{ a.s.},
\]
where $P_{n,k}$ is the $k$th largest $\mu$-measure of the collection of
all fringe subtrees of~$\mathcal{T}$ attached to $R_n$, or again
as
\[
S_{1/2} (P_n) = \lim_{m \rightarrow\infty}
m^{-1/2} K_{n,m},
\]
where $K_{n,m}$ is the number of junction vertices $J_i$ with $i \le m$ such
that $J_i \in R_n$.

\item[(iii)] Conditionally given $\lambda(R_n)$, the tree $R_n$
consists of
$2n-1$ segments, whose relative lengths, when listed in order of depth-first
search of $R_n$, passing first along $[\![0,\Sigma_1]\!] = R_1$, then along
$]\!]J_1,\Sigma_2]\!] = R_{2} \setminus R_1$, and so on, is independent
of $R_n$, with the same Dirichlet distribution with $2n-1$ parameters
equal to $1$ as the sequence of $2n-1$ spacings generated by a sequence
of $2n-2$ independent uniform variables on $[0,1]$.

\item[(iv)] For $n \ge2$ the combinatorial shape of $R_n$ is
equally likely to be
any of the $1 \times3 \times\cdots\times(2n-3)$ possible shapes of binary
trees with root $0$ and leaves labeled by $[n]$, independently of
$\lambda(R_n)$ and of the sequence of relative lengths of segments.

\item[(v)] Conditionally given $\lambda(R_n) = \ell_1$, the
combinatorial shape of
$R_n$ and the sequence of relative lengths of segments,
let $(\sigma_v, 0 \le v \le\ell_1)$ be a path which traverses $R_n$ at
unit speed, passing first along $[\![0,L_1]\!] = R_1$, then along
$]\!]J_1,L_2]\!] = R_{2} \setminus R_1$, and so on, and let $R_{n,\ell}$
be the range of
$(\sigma_v, 0 \le v \le\ell)$, so that by construction
$\lambda(R_{n,\ell}) = \ell$ for all $0 \le\ell\le\lambda(R_n) =
\ell_1$.
Then the cumulative mass process
$
(\mu_n ( R_{n,\ell} ), 0 \le\ell\le\lambda(R_n) \giv\lambda(R_n)
= \ell_1)
$
has the same distribution with exchangeable increments as the inverse of
the local time process at $0$ of a Brownian bridge $B^{\mathrm{br}}$
conditionally
on $L_1^0(B^{\mathrm{br}}) = \ell_1$.
\end{longlist}
\end{prop}

\begin{pf}
This can largely be read from the preceding discussion and known descriptions
of $(R_n,n\ge1)$ and properties of Poisson--Dirichlet distributions.
Partial results appear in many places, including
\cite{Ald-crt3}, Section~4.3, \cite{HMPW}, Proposition~18, and \cite
{PW09}, Proposition 14.

In the terminology of \cite{DGM}, part (i) is a particular case of
their result that crushing a size-biased pick from a ranked list
with distribution $\mathrm{PD}(\alpha,\theta)$ into $\mathrm
{PD}(\alpha,1-\alpha
)$-distributed proportions yields a
$\mathrm{PD}(\alpha,\theta+1)$ ranked vector. This can also be read from
Aldous's sequential description of the growth of $(R_n,n\ge1)$. That
description implies
part (v) quite easily. To deduce (ii) from (v), observe that $P_n$ is
the sequence of
ranked jumps of the cumulative mass process
$(\mu_n ( R_{n,\ell} ), 0 \le\ell\le\lambda(R_n))$ which
given $\lambda(R_n) = \ell$ is distributed like the ranked lengths of
excursion intervals
of a standard Brownian bridge. As $n$ changes, the distribution of
$P_n$ is therefore obtained
from that of $P_1$ by tilting the distribution by the density factor between
the exponential distribution of $\lambda(R_1)^2/2$ and the
Gamma$(n,1)$ distribution of
$\lambda(R_n)^2/2$. But this factor is just $(\lambda
(R_1)^2/2)^{n-1}/\Gamma(n)$ where
$\lambda(R_1)^2/2 = (S_{1/2}(P_1))^2/4$, which is precisely the
density factor between
$\mathrm{PD}(\frac{1}{2}, (2n-1)/2)$ and $\mathrm{PD}(\frac
{1}{2},\frac{1}{2})$; see \cite{csp}, Theorem 3.13.

As another check of the consistency of the two different descriptions
of $\lambda(R_n)$,
observe that the description in terms of independent exponential variables
gives
%
\begin{equation}
\label{rmoms} \mathbb{E}\bigl( \lambda(R_n)^\rho\bigr) =
2^{\rho/2 } \mathbb{E}\bigl( \Gamma_n^{\rho/2} \bigr) =
2^{p/2} \Gamma(n + \rho/2)/\Gamma(n),
\end{equation}
whereas for $S_{\alpha,\theta}$ the $\alpha$-diversity of
a $\mathrm{PD}(\alpha,\theta)$ random discrete distribution
it is known \cite{csp}
that
%
\begin{equation}
\label{mlmoms} \mathbb{E}\bigl( S_{\alpha,\theta}^ \rho\bigr) =
\frac{\Gamma( \theta/\alpha+ \rho+ 1 ) \Gamma(\theta+ 1 )}{
\Gamma(\theta+ p \alpha+ 1 ) \Gamma( \theta/\alpha+ 1 )}
\end{equation}
so in particular
%
\begin{equation}
\label{mlmomshf} \mathbb{E}\bigl( S_{1/2,(2n-1)/2}^ \rho\bigr) =
\frac{\Gamma( 2n + \rho) \Gamma(n + (1/2))} {\Gamma(n +
\rho/2 +(1/2)) \Gamma( 2n )} = \frac{ 2 ^\rho\Gamma(n + \rho/2 ) }{ \Gamma(n ) },
\end{equation}
where the second equality uses the Gamma duplication formula
$\Gamma(2z) = 2^{2z-1}\* \Gamma(z) \Gamma(z + \frac{1}{2}) /\Gamma
(\frac{1}{2})$.
\end{pf}\vspace*{-12pt}

%
\begin{appendix}\label{appA}%
\section*{Appendix: Proof of Lemma \texorpdfstring{\protect\ref{454654156489784}}{6}}
First note that the assumptions imply that $\tau$ is exponentially
distributed and, furthermore, that for all nonnegative Borel functions
$g$ and all $t\ge0$
\[
\mathbb{E}\bigl(g(F_\tau)1_{\{\tau>t\}}\bigr)=\mathbb{P}(\tau>t)
\mathbb{E}\bigl(g(F_{t+\tau-t})|\tau>t\bigr)=\mathbb{P}(\tau>t)\mathbb{E}
\bigl(g(F_\tau)\bigr),
\]
so that $\tau$ is independent of $F_\tau$. Since $\tau$ is a
stopping time, $(F_{\tau+s},s\ge0)$ is
independent of $\mathcal{F}_\tau$ and has the same distribution as
$(F_s,s\ge0)$. Consider a sequence
of independent copies $(F^{(n)}_s,0\le s\le\tau_n)$, $n\ge1$, of
$(F_s,0\le s\le\tau)$ and splice them together as
\[
(\widetilde{F}_{\tau_1+\cdots+\tau_{n-1}+s},\widetilde{m}_{\tau
_1+\cdots+\tau_{n-1}+s})=
\bigl(F^{(i)}_s,1_{\{s=\tau_n\}}\bigr),\qquad0<s\le
\tau_n,n\ge1.
\]
Then\vspace*{1.5pt} $((F_t,t\ge0),\tau)\stackrel{d}{=}((\widetilde{F}_t,t\ge
0),\widetilde
{\tau})$, by construction, where $\widetilde{\tau}=\tau_1$. Also by
construction,
$\widetilde{F}=(\widetilde{F}_t,t\ge0)$ and $\widetilde{F}{}^\bull
=((\widetilde{F}_t\widetilde{m}_t+(1-\widetilde{m}_t)),t\ge0)$ are
$\widetilde{\mathcal{F}}$-Poisson point processes in their joint
natural filtration $\widetilde{\mathcal{F}}$. The intensity measure of
$\widetilde{F}{}^\bull$ is the distribution of
$F_\tau$ times the rate of $\tau$. Since $F_\tau\neq1$ a.s., this
intensity measure is absolutely continuous with respect to $\Lambda$
as otherwise, we could find a Borel
set $A$ with $\mathbb{P}(F_\tau\in A)>0$ and $\Lambda(A)=0$, so the
restrictions of $\widetilde{F}$ and
$\widetilde{F}{}^\bull$ would reveal points of $\widetilde
{F}{}^\bull$ that are not points of $\widetilde{F}$, a contradiction.
We denote the Radon--Nikodym derivative by $K(u,dk)$ and set
$p(u)=K(u,\{1\})$, so that $\widetilde{F}{}^\bull$ is a Poisson
point process with intensity measure $p(u)\Lambda(du)$. We note that
$p(u)\le1$ for \mbox{$\Lambda$-a.e.} $u\in(0,1)$ as otherwise
restrictions to the Borel set $A=\{u\in(0,1)\dvtx  p(u)>1\}$ would
reveal points of $\widetilde{F}{}^\bull$ that are not points of
$\widetilde{F}$, a contradiction.

The proof is not complete yet because we have not yet shown that
$((\widetilde{F}_t,\widetilde{m}_t),\break t\ge0)$ is an
$\widetilde{\mathcal{F}}$-Poisson point process, or equivalently that
the unmarked points
$\widetilde{F}{}^\circ=((\widetilde{F}_t(1-\widetilde
{m}_t)+\widetilde{m}_t),t\ge0)$ form an $\widetilde{\mathcal
{F}}$-Poisson point
process. We can represent the point process $\widetilde{F}{}^\circ$ as
a random measure
\[
N^\circ_t(A)=\#\bigl\{s\le t\dvtx \widetilde{F}{}^\circ_s
\in A\bigr\},\qquad\mbox{$t\ge0$ and $A$ Borel subset of $(0,1)$}.
\]
For $A$ with $\Lambda(A)<\infty$ and $\int_Ap(u)\Lambda(du)>0$, we
now show that $(N_t^\circ(A),t\ge0)$ is an $\widetilde{\mathcal
{F}}$-Poisson process. Consider
$\overline{N}_t(A)=\#\{s\le t\dvtx \widetilde{F}_s\in A\}$, which we
know is an $\widetilde{\mathcal{F}}$-Poisson process with rate
$\Lambda(A)$. The time $\tau^A=\inf\{t\ge0\dvtx \widetilde
{F}{}^\bull_t\in A\}$ is an
\mbox{$\widetilde{\mathcal{F}}$-}stopping time, since $\widetilde
{F}{}^\bull$ is an $\widetilde{\mathcal{F}}$-Poisson point process.
Denote by
$(T_n,n\ge1)$ the times of the points of $\overline{N}(A)$, also
$\widetilde{\mathcal{F}}$-stopping times. Let $A_n=\{\tau^A=T_n\}$
and set $q=\mathbb{P}(\tau^A=T_1)$. Then, by
the strong Markov property of $\overline{N}(A)$ at $T_n$, we find that
also for $n\ge2$
%
\begin{equation}
\label{thinn}\qquad \mathbb{P}(A_n)=\mathbb{P}\bigl(\tau^A\neq
T_1\bigr)\mathbb{P}\bigl(A_n|\tau^A>T_1
\bigr)=(1-q)\mathbb{P}(A_{n-1})=(1-q)^{n-1}q.
\end{equation}
Hence,\vspace*{2pt} we find that $\tau^A=T_G$ for a $G$ geometric with parameter
$q$, which is a stopping time in the discrete filtration
$(\widetilde{\mathcal{F}}_{T_n},n\ge1)$. By Wald's equation,
\[
\mathbb{E}\bigl(\tau^A\bigr)=\mathbb{E}(G)\mathbb{E}(T_1)
\qquad\mbox{hence }\int_A p(u)\Lambda(du)=q\Lambda(A).
\]
From (\ref{thinn}) and the strong Markov property at $\tau^A$, we deduce
that $(N^\bull(A),\break \overline{N}(A))$ are such that
$N^\bull(A)$ is a $q$-thinning of $\overline{N}(A)$, so $N^\circ
(A)=\overline{N}(A)-N^\bull(A)$ is also a Poisson process
with rate
\[
(1-q)\Lambda(du)=\int_A \bigl(1-p(u)\bigr)\Lambda(du).
\]
By \cite{kal}, Theorem 12.8, this suffices to identify $N^\circ$ as a
Poisson random measure with intensity measure
$(1-p(u))\Lambda(du)$.
\end{appendix}



%

\printaddresses
\end{document}